\RequirePackage{fix-cm}
\documentclass[smallextended]{svjour3}       
\smartqed  
\usepackage{graphicx}
\emergencystretch=\maxdimen
\usepackage{epsfig}
\usepackage{graphicx}
\usepackage{makeidx}
\usepackage{pgfplots}

\usepackage{mwe}
\usetikzlibrary{quotes,angles}
\usepackage{epsfig}
\usepackage{latexsym}
\usepackage{graphicx}
\usepackage{makeidx}

\usepackage{amsmath}
\usepackage{amssymb}

\usepackage{graphicx}
\usepackage{color}
\usepackage{epsfig}
\usepackage{lineno}
\usepackage{pict2e}
\usepackage[percent]{overpic}

\usepackage{multicol}
\usepackage{mwe}
\usepackage{subfigure}
\usepackage{graphicx}
\usepackage{makeidx}
\usetikzlibrary{quotes,angles}
\usepackage{latexsym}
\usepackage{graphicx}
\usepackage{amsmath}
\usepackage{amsthm}
\usepackage{color}
\usepackage{lineno,hyperref}
\usepackage{multirow}
\usepackage{rotating}
\usepackage{mathtools}

\usepackage{multicol}
\newcommand{\RN}[1]{%
  \textup{\uppercase\expandafter{\romannumeral#1}}%
}

\newtheorem{thm}{Theorem}
\newtheorem{lem}{Lemma}

\pgfplotsset{compat=1.17}

\begin{document}

\title{  Additive Average Schwarz Method for Elliptic Mortar Finite Element Problems with Highly Heterogeneous Coefficients 
}


\author{Ali Khademi \and Leszek Marcinkowski \and Sanjib Kumar Acharya \and Talal Rahman 
}


\institute{Ali Khademi \at
              Department of Computer science, Electrical engineering and Mathematical sciences, \\
       Western Norway University of Applied Sciences, P.O. Box 7030, Bergen, Norway \\
              \email{ali.khademi@hvl.no, akhademi.math@gmail.com}           
           \and
           Leszek Marcinkowski \at
              Department of Mathematics, Warsaw University, Poland\\ \email{lmarcin@mimuw.edu.pl}
              \and
              Sanjib Kumar Acharya \at
              Institute of Chemical Technology Mumbai, Indian Oil Campus, Odisha, Bhubaneswar, India\\
             \email{acharya.k.sanjib@gmail.com}
             \and 
             Talal Rahman \at
             Department of Computer science, Electrical engineering and Mathematical sciences, \\
       Western Norway University of Applied Sciences, P.O. Box 7030, Bergen, Norway \\
             \email{talal.rahman@hvl.no}
}

\date{Received: date / Accepted: date}

\maketitle

\begin{abstract}
 In this paper, we extend the additive average Schwarz method to solve second order elliptic boundary value problems with heterogeneous coefficients inside the subdomains and across their interfaces by the mortar technique, where the mortar finite element discretization is on nonmatching meshes. In this two-level method, we enrich the coarse space in two different ways, i.e., by adding eigenfunctions of two variants of the generalized eigenvalue problems. We prove that the condition numbers of the systems of algebraic equations resulting from the extended additive average Schwarz method, corresponding to both coarse spaces, are of the order  $O(H/h)$ and independent of jumps in the coefficients, where $H$ and $h$ are the mesh parameters.

\keywords{Domain decomposition method \and Additive Schwarz method \and Generalized eigenvalue problem \and Mortar finite elements }
\end{abstract}

\section{Introduction}

Domain decomposition methods are efficient and powerful iterative methods to solve large algebraic systems arising from a finite element discretization of elliptic boundary value problems \cite{Smith-book-1996,Tallec-1991,Xu-1998}. They can also be regarded as a procedure of producing preconditioners for other iterative methods, such as the conjugate gradient method, for achieving fast convergence. In both approaches, to solve an original problem defined on a bounded Lipschitz domain $\overline{\Omega}= \cup_{i=1}^{N}\overline{\Omega}_i$, it is equivalent to solve many subproblems defined locally on the subdomains $\Omega_i$ in parallel. To obtain fast convergence, Dryja and Widlund \cite{Dryja-1987}, and Matsonki and Nepomnyaschikh \cite{Matsokin-1985} proposed to add one global problem and introduced the additive Schwarz methods to solve the global and local problems in parallel. 

Corresponding to the global and local problems, different coarse and fine spaces can be constructed, for instance,  see \cite{Brenner-1996,Bjorstad-1996,Cowsar-1993,Sarkis-1993}. It is worth mentioning that coarse spaces are more important than fine spaces since they have a key role in the central of scalability for domain decomposition methods. Therefore, in \cite{Bjorstad-1996} one of the simplest and efficient ways to construct coarse space, called the average coarse space, was proposed. Consequently, the two-level additive average Schwarz method was introduced and developed to solve many kinds of elliptic problems with the continuous and discontinuous coefficients, see \cite{Bjorstad-2003,Dryja-2010,Feng-2002,Loneland-2016}. Hence, an extension of the additive average Schwarz method to solve elliptic model problems arising from many applications such as composite materials with highly heterogeneous coefficients is of particular interest in this paper because it has a high-level performance.

 In general, in terms of distributions of the coefficients, we can classify elliptic problems with the heterogeneous coefficients into two classes, i.e., elliptic problems with jumps in the coefficients only inside the subdomains and ones with jumps in the coefficients both inside the subdomains and on the subdomain interfaces. For the first class, where the heterogeneous coefficients are piecewise constants for subdomains $\Omega_i, i=1,\ldots,N$, the additive Schwarz method was developed and analyzed in \cite{Bjorstad-1996,Chan-1994,Sarkis-1997,Toselli-2005} and references therein. For the second class with the large variation in the heterogeneous coefficients, the classical coarse spaces lead the condition numbers of the preconditioned systems to blow up.  Consequently, the convergence rates of iterative methods will deteriorate. To alleviate this difficulty, the coarse spaces can be enriched by combining their structure with spectral spaces.

 The idea of the coarse spectral space was introduced in \cite{Brezina-1999} and extended as the spectral algebraic multigrid method in \cite{Efandiev-2011}. The spectral construction of this new space is achieved by solving the generalized eigenvalue problems locally. Due to its crucial role in preventing the impact of large jumps in the coefficients on bounds of the condition numbers of preconditioners, the new coarse space has received considerable attention.  As a result, it is extensively developed for overlapping Schwarz methods \cite{Eikeland-2019,Gander-2015,Heinlein-2019,Spillane-2014}, nonoverlapping additive Schwarz method \cite{Marcinkowski-2018}, balancing domain decomposition methods \cite{Kim-2015,Kim-2017,Mandel-2012,Spillane-2013}, and  nonlinear domain decomposition \cite{Klawonn-2014}, where all  nodes on the subdomains interfaces for all references are matching grids. Hence, in this paper, we focus on the nonmatching grids, which are common due to heterogeneous materials in real-life problems. In practice, one may use different triangulations for polygonal subdomains independent of other triangulated subdomains. More precisely, there exist situations, where two subdomains with a common interface have fine and coarse (or fine but with different mesh sizes compared to other) triangulations. Therefore, the nonmatching grids on the subdomains interfaces are unavoidable and cause  consistency errors in numerical methods, which can be handled by using the mortar techniques \cite{Arbogast-2000,Arbogast-2007,Arbogast-2013,Bernardi-2005,Seshiyer-2000}. Hence, the additive Schwarz method for the mortar finite element method was introduced  \cite{Bjorstad-2003}, and modified \cite{Rahman-2005} with Crouzeix-Raviart mortar finite elements.
 
 The main aim of this paper is to enrich the coarse space used in the classical nonoverlapping additive average Schwarz method by using the idea of coarse spectral space for elliptic problems with highly heterogeneous coefficients inside the subdomains and across their interfaces where the mortar finite element discretization is on the nonmatching meshes. To achieve this, our new coarse space consists of two subspaces. The first subspace is the common  coarse space in the classical additive average Schwarz method used in \cite{Bjorstad-2003}, i.e., for the fixed $i$, is the range of a linear operator defined on $\Omega_i$  such that it is either the nodal values of a function $u \in V_h$ inside $\Omega_i$ or the average of nodal values of the function $u$ on the mortar and nonmortar sides of $\Omega_i$, where $V_h$ is a finite space of P1 conforming elements defined on a fine triangulation of $\Omega$ and vanishing on $\partial \Omega$. The second subspace has a particular spectral structure. To obtain the basis functions for this subspace, we solve the generalized eigenvalue problems restricted to each subdomain as in \cite{Marcinkowski-2018}.
 
 To define the proper generalized eigenvalues problems, we require determining minimum values of the coefficients over each subdomain's triangulation to estimate the condition number of additive average Schwarz preconditioners independent of the large eigenvalues caused by the large jumps in the coefficients. Motivated by the ideas from \cite{Marcinkowski-2018}, we consider two different types of the generalized eigenvalue problems based on either minimum of the coefficients over the whole subdomain or just minimum of ones over the layer connected to the boundary of the subdomain with one vertex or with one edge of the triangles inside of the subdomain. Solving these generalized eigenvalue problems lead to finding orthogonal basis functions enriching the coarse spaces. With these new coarse spaces, we prove that the condition numbers of the produced preconditioners are of the order $O(H/h)$ and independent of the number of subdomains. Due to the definition of the second type layer, it has faster performance than the first one through the implementations with numerical-software packages. See also Section \ref{SNumericalResults} for numerical results.

  
 The outline of this paper is as follows. In Section \ref{SDisProb}, after introducing a discrete problem, we define the mortar condition and the space of basis functions satisfying that condition. Furthermore, several figures related to those functions are also given there. Section \ref{AddAveSchw} is devoted to introducing the additive average Schwarz method, where the average interpolation operator has two different types. This operator consists of the natural extension of the standard average interpolation operator for mortar case and orthogonal operators defined in the next section. In Section \ref{EnrichCorse}, the generalized eigenvalue problems in terms of how the minimum values of jumps in the coefficients over the subdomains' triangulations can be defined are introduced. In Section \ref{CondNumb}, following the standard additive Schwarz framework \cite{Smith-book-1996}, we drive an optimal estimate of the condition number of the produced preconditioners with the aid of removing bad eigenvalues, which are influenced by the large jumps in the coefficients. Finally, to verify our theoretical results' validity, in Section \ref{SNumericalResults} some numerical experiments are reported.
 
\section{Discrete problem}\label{SDisProb}

Let $\Omega \subset \mathbb{R}^2$ be a bounded Lipschitz domain with a nonoverlapping partition $\{ \Omega_i \}^{N}_{i=1}$ of polygonal  subdomains such that $\overline{ \Omega} = \cup^{N}_{i=1} \overline{\Omega}_i$.
We consider an elliptic model problem defined on $\Omega$: Find $u_*\in H^{1}_{0}(\Omega)$ such that
\begin{equation}\label{eq:general_model}
a(u_*,v)=f(v),~v\in H^1_0(\Omega),
\end{equation}
where $$a(u,v)=\sum^{N}_{i=1}a_i(u,v)=\sum^{N}_{i=1}(\alpha_i \nabla u,\nabla v)_{L^{2}(\Omega_i)}$$
and
$$f(v)=\int_{\Omega}fv~dx=\sum^{N}_{i=1}\int_{\Omega_i}fv ~dx.$$

Here, $f\in L^{2}(\Omega)$, $\alpha\in L^{\infty}(\Omega)$ and $\alpha_i(x)$ is the restriction of $\alpha(x)$ over $\Omega_i$. Further, we assume that there exists a positive constant $\alpha_0$ such that $\alpha(x)>\alpha_0$. 

The partition $\{\Omega_i\}^{N}_{i=1}$ forms a coarse triangulation  of $\Omega$ with the mesh parameter  $H=\max\{  H_i, ~i=1,\ldots,N\}$, where $H_i$ are diameters of $\Omega_i$. We assume this partition to be geometrically conforming,  i.e.,  $\partial\Omega_{i}\cap\partial\Omega_{j}$
$(i\neq j)$ is  a vertex  or a whole edge of both subdomains $\Omega_{i}$ and $\Omega_{j}$ or is empty. Further, we denote the set of all vertices of such coarse triangulation, except those belonging to $\partial \Omega$, by $\mathcal{N}_H$. We also denote the triangulation of the subdomain 
\begin{figure}[!ht]
\setlength{\unitlength}{0.1\textwidth}
\hspace{-3.2cm}
\begin{picture}(10,4.5)
\put(0,0){ \includegraphics[width=17cm]{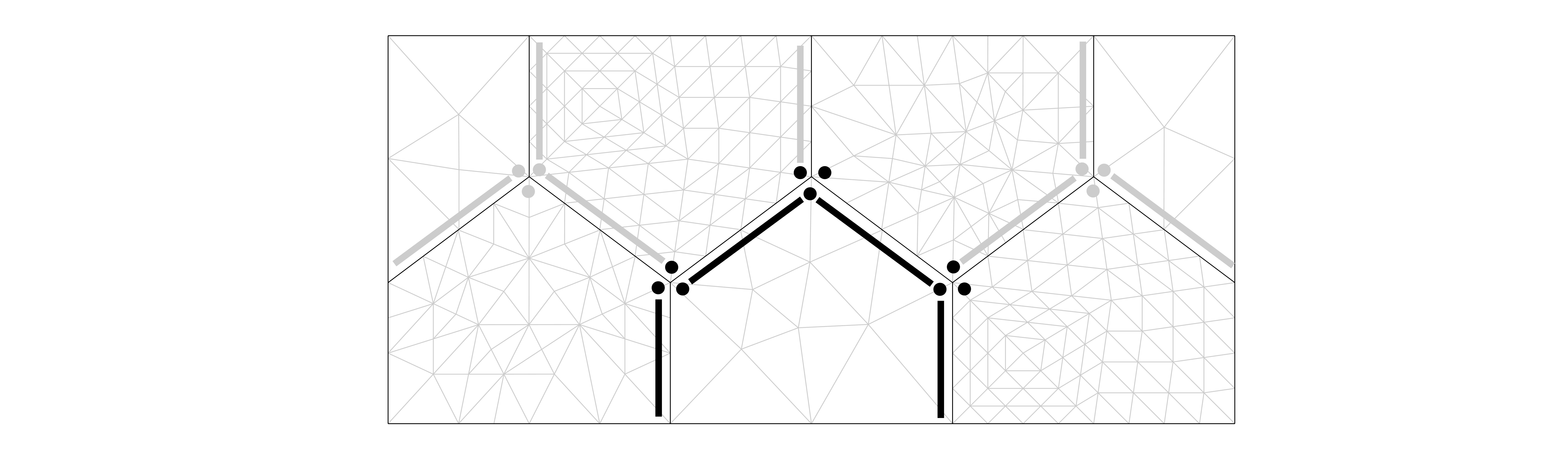}}
  \put(6.6,0.7){$\Omega_i$}
  \put(9.2,3.5){$\Omega_j$}
  \put(7.5,1.9){$\gamma_{m(i)}$}
  \put(8,2.4){$\delta_{m(j)}$}
  \end{picture}
  \caption{ The mortar and nonmortar sides of subdomains $\Omega_i$ and $\Omega_j$ which are  denoted by  $\gamma_{m(i)}$ and $\delta_{m(j)}$, respectively.}
  \label{fig: mortar and nonmortar}
\end{figure}
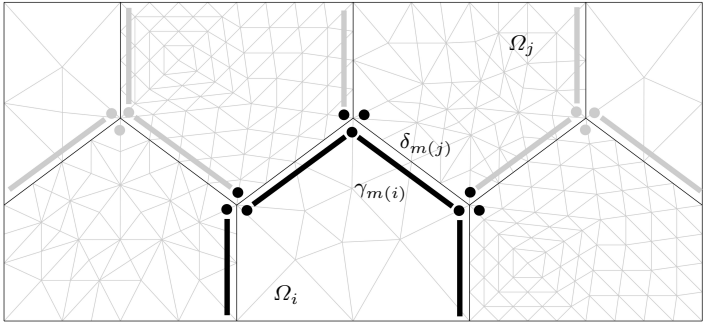
$\Omega_i$ by $\mathcal{T}_{i}$ which consists of triangles satisfying the shape regular property \cite{Ciarlet-1978} inside of $\Omega_i$, and quasi-uniform triangles touching $\partial \Omega_i$ with the mesh size $h_i$. Further, we denote  the set of all internal nodes of the fine triangulation $\mathcal{T}_i$ by $\imath_i$. 
 
We define the product space $X_h$ on the computational domain $\Omega$ by
$$X_h(\Omega)=X_1(\Omega_1) \times \ldots \times  X_N(\Omega_N),$$
where $X_i(\Omega_i), i=1,\ldots,N$ are the finite element spaces of the continuous piecewise linear functions defined on $\mathcal{T}_i$ and vanishing on $\partial\Omega\cap\partial\Omega_i$.  We denote all nodal points on $\bigcup_{i=1}^{N} \mathcal{T}_i$ except those on $\partial \Omega$ by $\mathcal{N}_{\Omega}$. Further, we denote  the set of basis functions associated with the set of nodal points $\mathcal{N}_{\Omega}$ by $\{\phi_l \}_{l \in \mathcal{N}_{\Omega}}$.
 
 Due to independent triangulations inside of each subdomain, on each side of the interface $\overline{\Gamma}_{ij}=\overline{\Omega}_i \cap \overline{\Omega}_j$ we may have different discretization (cf. Figure~\ref{fig: mortar and nonmortar}). We select one side of $\Gamma_{ij}$ as the mortar side, and the other side as the nonmortar side denoted by $\gamma_{m(i)}$ and $\delta_{m(j)}$, respectively. It is obvious that $\Gamma_{ij}=\gamma_{m(i)}=\delta_{m(j)}$. Further, we denote the nodes on the mortar and nonmortar sides by $m_0,m_1,\ldots,m_{n_m+1}$ and $s_0,s_1,\ldots,s_{n_s+1}$, respectively. 
 
Let $W^{h_i}(\Gamma_{ij})$ and $W^{h_j}(\Gamma_{ij})$ denote the restrictions of $X_i(\Omega_i)$ and $X_{j}(\Omega_j)$ onto $\Gamma_{ij}$, respectively. Now, the nonmatching grids on the subdomain interfaces impose discontinuities for the functions belong to $X_h$. Therefore, we need to define a week continuity condition. To this end, we first define the  projection $\Pi_m(u_i, \mathrm{Tr}~ v_j) : L^{2}(\delta_{m(j)})\rightarrow W^{h_j}(\delta_{m(j)})$ by
\begin{equation}\label{mortcond}
\int_{\delta_{m(j)}}
\Big( u_{i_{|_{\gamma_{m(i)}}}} - ~
\Pi_m(u_i, \mathrm{Tr}~ v_j) \Big) \psi ~ dx=0,  \quad \gamma_{m(i)}=\Gamma_{ij}=\delta_{m(j)},
\end{equation}
for all functions $\psi\in M^{h_j}(\delta_{m(j)})) = \mathrm{span} \{ \varphi_l \}_{l=0}^{n_s+1}$ and 
\begin{equation*}
 \Pi_m(u_i, \mathrm{Tr}~ v_j)_{|_{\partial \delta_{m(j)}}}= v_{j_{ |_{\partial \delta_{m(j)}}}}, \end{equation*}
where $u_i$, $\mathrm{Tr}~ v_j$, and $M^{h_j}(\delta_{m(j)})$ are the restriction of $u$ into each subdomains $\Omega_i, i=1,\ldots,N$, the trace of $v_j$, and a subspace of $W^{h_j}(\delta_{m(j)})$ with constant values on the elements touching $\partial\delta_{m(j)}$, more precisely, at two end points of $\delta_{m(j)}$, respectively. Now, we say a function $u_h =\{u_i\}^{N}_{i=1}\in X_{h}$ satisfies the  
\begin{figure}[!ht]
\centering
\setlength{\unitlength}{0.1\textwidth}
\vspace{6cm}
\begin{picture}(15,6.5)
\put(0,7.5){\includegraphics[width=5.8cm]{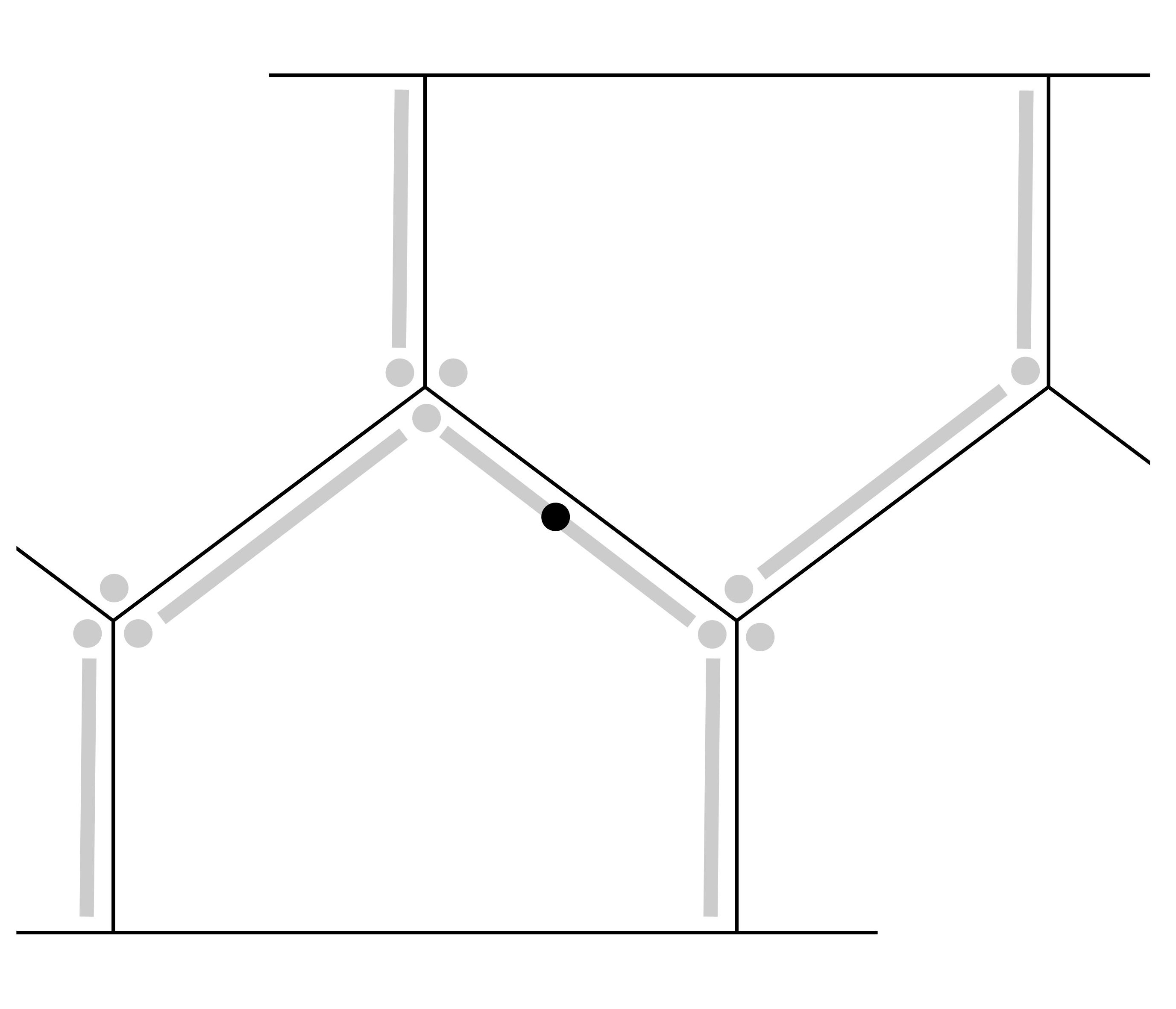}}
  \put(.6,8.1){ $\Omega_i$}
  \put(3.6,11){ $\Omega_j$}
  \put(1.8,9.5){ $x_k$}
  \put(2.2,9.9){$\Pi_{m}$}
  \put(2.2,7.3){$(\RN{1})$}
  \put(5,7.5){\includegraphics[width=5.8cm]{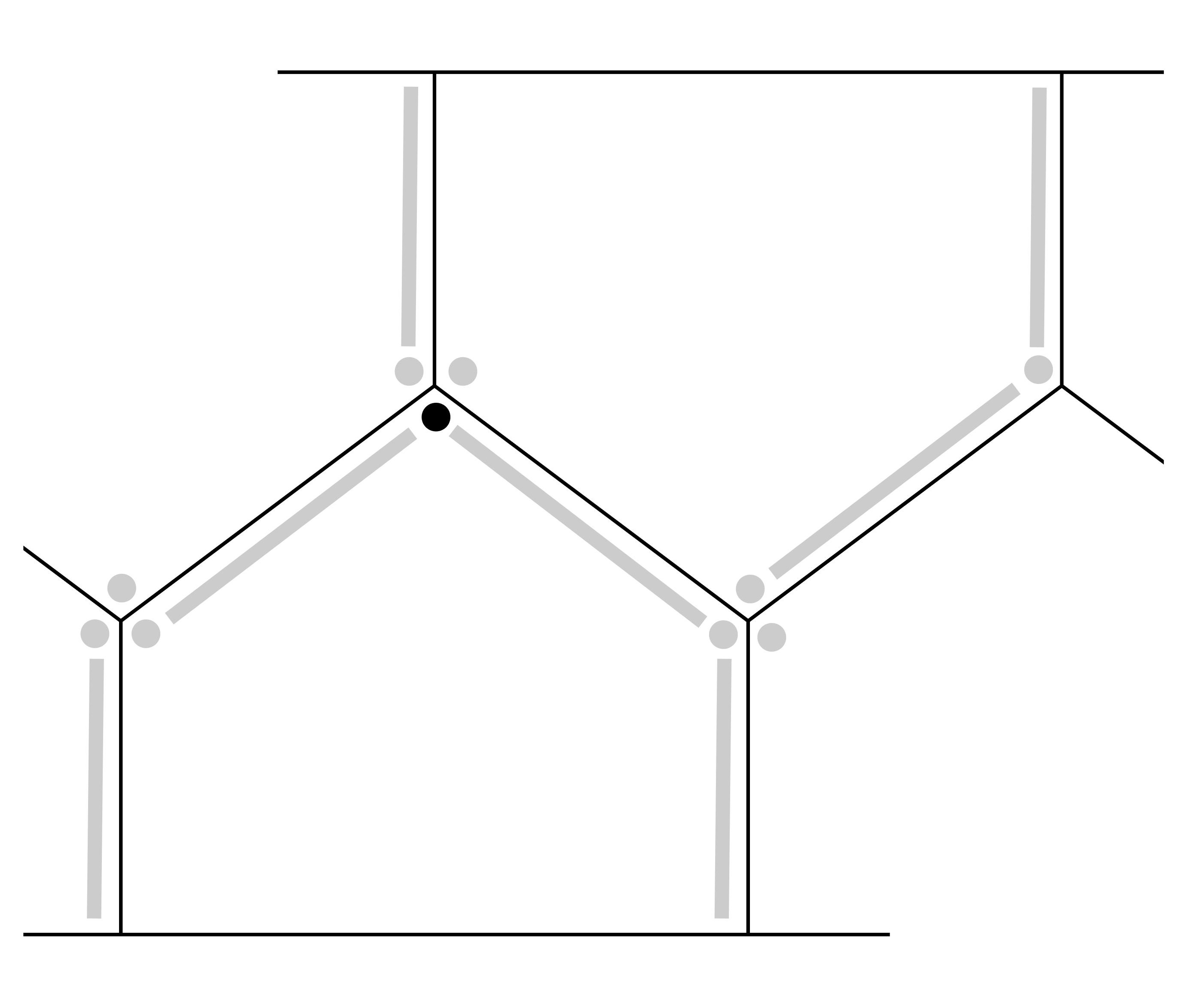}}
  \put(5.6,8.1){ $\Omega_i$}
  \put(8.6,11){ $\Omega_j$}
  \put(6.55,9.6){ $x_k$}
  \put(6.,9.9){$\Pi_{n}$}
  \put(7.3,9.9){$\Pi_{m}$}
  \put(6.,9.1){$\gamma_{n(i)}$}
  \put(7,9.1){$\gamma_{m(i)}$}
  \put(7.2,7.3){$(\RN{2})$}
\put(0,2.5){\includegraphics[width=5.8cm]{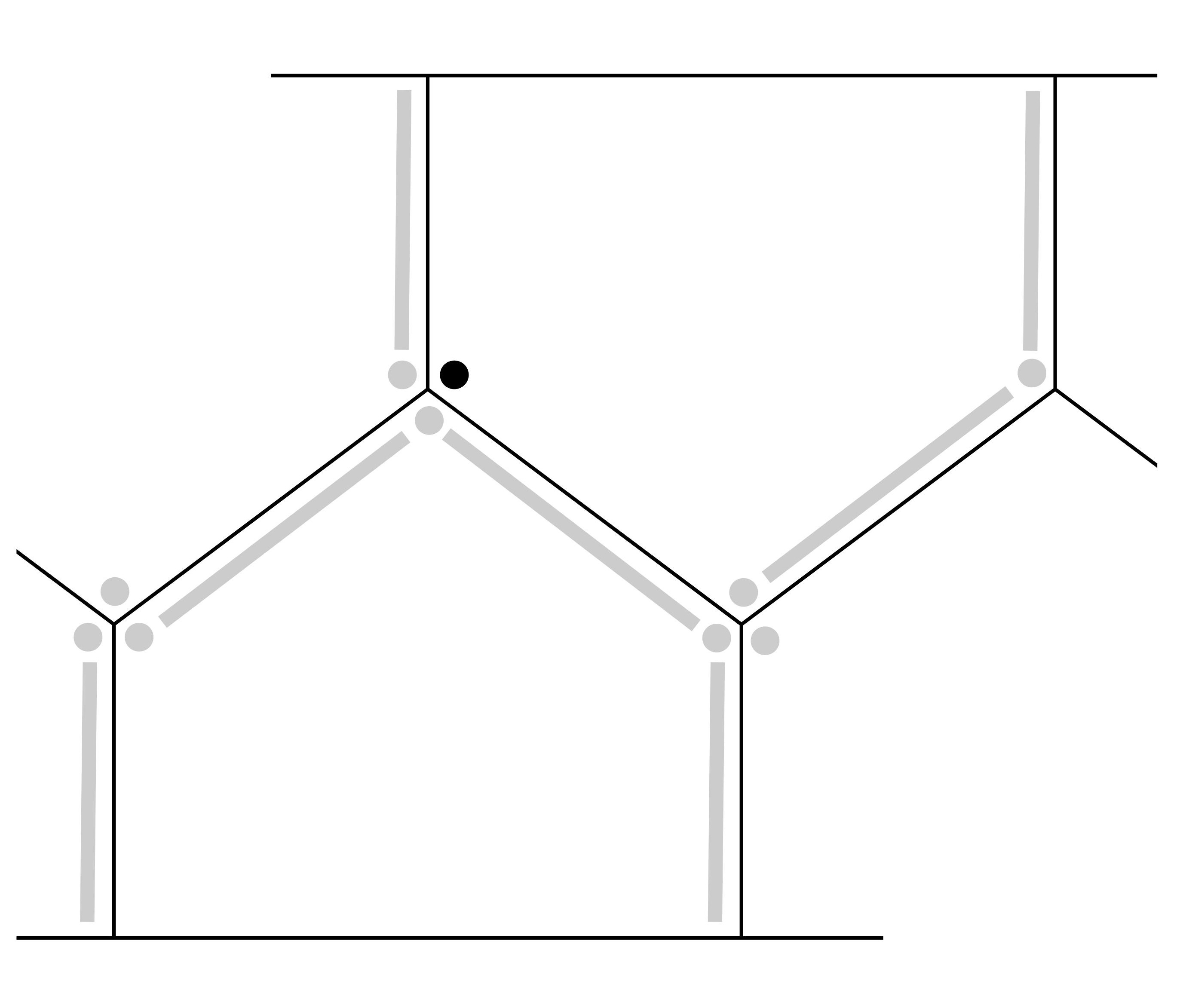}}
\put(.6,3.05){ $\Omega_i$}
\put(3.6,6){ $\Omega_j$}
\put(1.9,5.8){$\Pi_{n}$}
\put(1.,5.5){$\gamma_{n(p)}$}
\put(1.95,5.15){ $x_k$}
\put(1.75,4.25){ $\gamma_{m(i)}$}
\put(2.4,4.8){$\Pi_{m}$}
\put(2.2,2.3){$(\RN{3})$}
  \put(5,2.5){ \includegraphics[width=5.8cm]{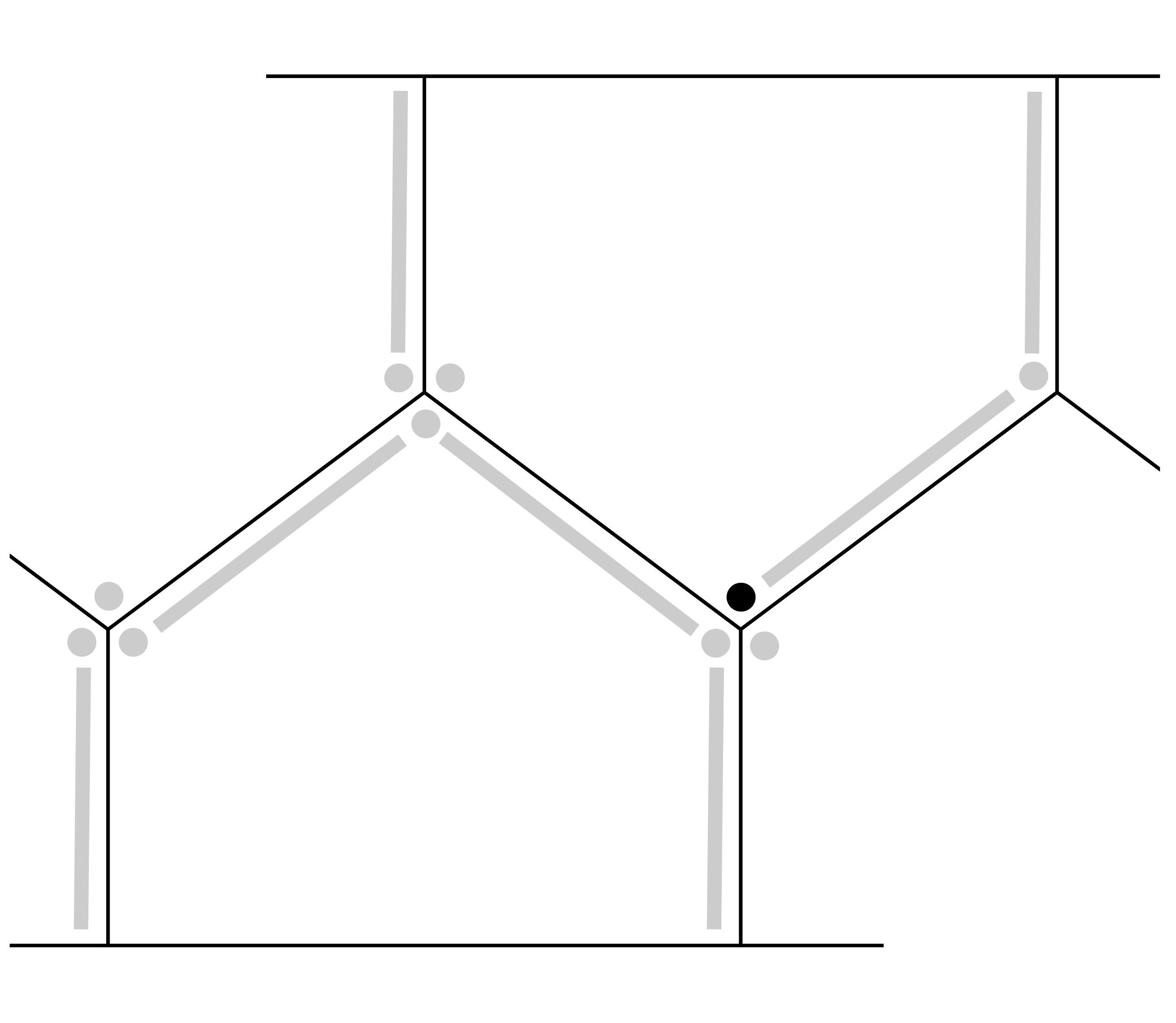}}
  \put(5.7,3.05){ $\Omega_i$}
  \put(8.7,6){ $\Omega_j$}
  \put(7.8,4.5){ $x_k$}
  \put(7.2,2.3){$(\RN{4})$}
  \put(7.,4.25){ $\gamma_{m(i)}$}
  \put(7.2,4.9){$\Pi_{m}$}
  \put(8.2,4.9){$\gamma_{n(j)}$}
  \put(8.7,4.2){$\Pi_{n}$}
  \end{picture}
  \vspace{-3.2cm}
  \centering
  \caption{The node $x_k$ is represented by the black and thick dot, where the value of the basis function $\phi_k(x)$ associated with $x_k$ is 1. Figures (\RN{1})-(\RN{4}) illustrate the different positions of $x_k$.}
  \label{Disterbution of node or vertex}
\end{figure}
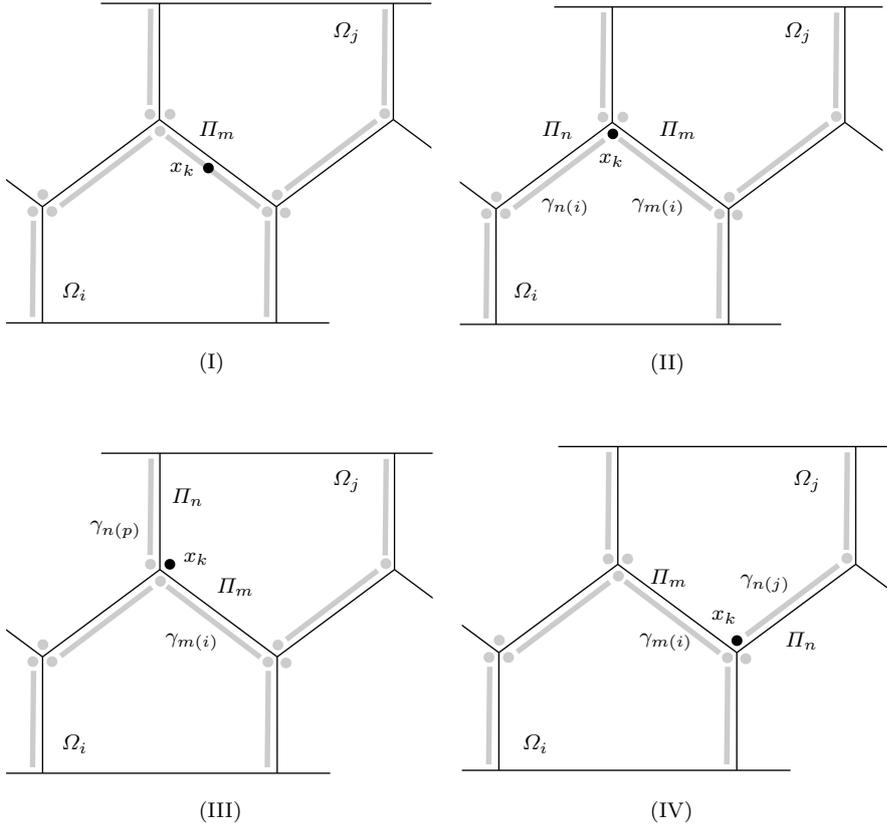
mortar condition on $\delta_{m(j)}$, if (\ref{mortcond}) holds. 
We denote the mortar space by $V_h$ in terms of the mortar condition, i.e.,
$$
V_h=\{  u_h,v_h \in X_h ~|~\Pi_m(u_h, \mathrm{Tr}~ v_h)=0, ~ \forall ~  \Gamma_{ij},~ i,j=1,\ldots,N \}.
$$

To see the structure of the basis functions spanning the mortar space $V_h$, let $\phi_k^{(i)}(x)$  be a basis function defined on $\Omega_i$ and  define $V_h=\mathrm{span}\{\phi_k\}$, where $\phi_k(x)$ is a basis function associated with a node $x_k$. Due to different positions of $x_k$, $\phi_k(x)$ takes different forms as follows (cf. Figure~\ref{Disterbution of node or vertex}).

\begin{enumerate}
    \item $x_k \in \imath_i:$   
    \hspace{-3cm}{$$\hspace{-7.5cm} \phi_k(x)=\phi_{k}^{(i)}(x).$$}

    \item $x_k \in \{ m_1,\ldots,m_{n_m} \}:$ 
    $$\phi_k(x)=\begin{cases}
    \phi_k^{(i)}(x), \hspace{1.85cm} \text{on } \overline{\Omega}_i,\\ 
    \Pi_m(\phi_k^{(i)}(x),0)(x), \hspace{.23cm} \text{on } \overline{\delta}_{m(j)}, \text{ where } \gamma_{m(i)}=\delta_{m(j)},\\
    0, \hspace{2.7cm} \text{otherwise.} 
    \end{cases}$$
    \item  $x_k \in \mathcal{N}_H:$
    \begin{enumerate}
        \item $x_k$ is a point of intersection between two mortar sides $\gamma_{m(i)}$ and $\gamma_{n(i)}$:
        $$\hspace{-2.7cm}\phi_{k}(x)= \begin{cases}
        \phi_k^{(i)}(x), \hspace{1.85cm} \text{on } \overline{\gamma}_{m(i)} \text{ and } \overline{\gamma}_{n(i)},\\ 
        \Pi_m(\phi_k^{(i)},0)(x), \hspace{0.7cm} \text{on } \overline{\delta}_{m(j)}, \\
        \Pi_n(\phi_k^{(i)},0)(x), \hspace{0.8cm}  \text{on } \overline{\delta}_{n(j)},\\
        0, \hspace{2.7cm} \text{otherwise.}
        \end{cases}$$
        
        \item $x_k$ is a point of intersection between mortar side $\gamma_{m(i)}$ and nonmortar side $\delta_{n(i)}$:
        $$\hspace{-4.3cm}\phi_{k}(x)=\begin{cases}
        \phi_k^{(i)}(x), \hspace{1.9cm} \text{on } \overline{\gamma}_{m(i)},\\ 
        \Pi_m(\phi_k^{(i)},0)(x), \hspace{0.65cm} \text{ on } \overline{\delta}_{m(j)},\\
        \Pi_n(0, \mathrm{Tr}~ \phi_k^{(i)})(x), \hspace{0.4cm} \text{on  }\overline{\delta}_{n(i)}, \\
        0, \hspace{2.8cm} \text{otherwise.}
        \end{cases}$$
        \item $x_k$ is a point of intersection between two nonmortar sides $\delta_{m(i)}$ and $\delta_{n(i)}$:
        $$\hspace{-2.7cm}\phi_{k}(x)=\begin{cases}
        \phi_k^{(i)}(x), \hspace{2cm} \text{on } \overline{\gamma}_{m(j)} \text { and } \overline{\gamma}_{n(j)},\\
        \Pi_m(0, \mathrm{Tr} ~ \phi_k^{(i)})(x), \hspace{0.4cm} \text{on } \overline{\delta}_{m(i)},\\
         \Pi_n(0, \mathrm{Tr} ~ \phi_k^{(i)})(x), \hspace{0.5cm} \text{on } \overline{\delta}_{n(i)},\\
         0, \hspace{2.9cm} \text{otherwise.}
        \end{cases}$$
        \end{enumerate}
\end{enumerate}

We now express the main problem as in the following form: Find $u^{*}_h = \{ u^h_i \}_{i=1}^{N} \in V_h$ such that
\begin{equation}\label{MainLinearForm}
a_h(u^{*}_h,v_h)=f(v_h) \quad \forall v_h \in V_h.
\end{equation}
  
In what follows, we consider the following matrix representation of the linear systems arising from the discretization of \eqref{MainLinearForm}.

$$\mathbf{A} \mathbf{v} = \mathbf{f}, $$
where $\mathbf{v}$ is the vector of all unknown  coefficients defined on $\overline{\Omega}$. Further, we consider the submatrices $\mathbf{A}_{\Omega_i} = \mathbf{R}_i\mathbf{A}\mathbf{R}_i^T$, where  $\mathbf{R}_i, i=1,\ldots,N$ are the restriction matrices such that $\mathbf{v}_i= \mathbf{R}_i \mathbf{v}$ are the vectors of coefficients defined on $\Omega_i \setminus \partial \Omega_i$. 
 
Employing the mortar condition, we  can compute some coefficients of $\mathbf{v}$ in terms of other ones. More precisely, where $\Gamma_{ij}=\gamma_{m(i)}=\delta_{m(j)}$ assume $\nu_s=(v_h(s_i))_{i=1}^{n_s}, \nu_c=(v_h(s_0),v_h(s_{n_s +1}))^T,\nu_m=(v_h(m_i))_{i=0}^{n_m+1}$ and consider the following matrix representations as in \cite{Dryja-2004}.
\begin{align*}
    \mathbf{M}_{\gamma} & :=\left( (\varphi_{s_i} , \phi_{m_j})_{L^2(\gamma)} \right) _{i=1;j=0}^{n_s ~~n_m+1} \in \mathbb{R}^{n_s \times(n_m+2)},\\
   \mathbf{S}_{\gamma} & := \left( (\varphi_{s_i} , \phi_{s_j})_{L^2(\gamma)} \right) _{i=1;j=1}^{n_s ~~n_s} \in \mathbb{R}^{n_s \times n_s},\\ 
   \mathbf{C}_{\gamma} & := \left( (\varphi_{s_i} , \phi_{m_j})_{L^2(\gamma)} \right) _{i=1;j=0,n_s+1}^{n_s} \in \mathbb{R}^{n_s \times2}.
\end{align*}

Hence
$$
\nu_s =\mathbf{S}_{\gamma}^{-1} (\mathbf{M}_{\gamma}\nu_m - \mathbf{C}_{\gamma}\nu_c).
$$

Consequently, in the next sections, we will focus only on all mortar, corner, and interior nodes. 
 
\section{The additive average Schwarz method}\label{AddAveSchw}

Let $V^{(i)}$, $i=1,\ldots,N$ be decompositions of $V_h$ that are restrictions of $V_h$ to $\Omega_i$ with zero on $\partial\Omega_i$ and on the subdomains $\Omega_j,~ j\neq i$. Further, let $V^{(type)}_{0}$ stands for two different types of the coarse spaces, distinguished by two notations $\RN{1}$ and $\RN{2}$, such that
$$V_h=V^{(type)}_{0}+\sum^{N}_{i=1}V^{(i)}, \quad type\in\{\RN{1},\RN{2}\},$$
where
$$V_0^{(type)}=I_0V_h+\sum^{N}_{i=1}W^{type}_{i},~~type\in\{\RN{1},\RN{2}\}.$$ 

We define the operator $I_0$ and the spaces $W^{type}_{i},~~type\in\{\RN{1},\RN{2}\}$ in the next section. Here, we use the exact bilinear form for all our subproblems. 
Thus, we have
 $$a(u^h,v^h)=a_i(u^h_i,v^h_i), \quad i=1,\ldots,N,$$
where $u^h=\{u^h_i\}^{N}_{i=1}\in V^{(i)}$ and $v^h=\{v^h_i\}^{N}_{i=1}\in V^{(i)}$. We define the projection-like operators $T^{(i)}:V_h\rightarrow V^{(i)}$ $i=1,\ldots,N$ and  $T^{(type)}_{0}:V_h\rightarrow V^{(type)}_{0}$ such that for $u^h\in V_h$, $T^{(i)}u^h$ and $T^{(type)}_{0} u^h$ are the solutions of
$$\hspace{0.3cm} a(T^{(i)}u^h,v^h)=a(u^h,v^h), ~~v^h\in V^{(i)}, ~~ i=1,\ldots,N$$
and
\begin{equation*}
a(T^{(type)}_{0}u^h,v^h)=a(u^h,v^h), \quad v^h \in V_{0}^{(type)}, \quad type \in \{ \RN{1},\RN{2} \}.
\end{equation*}

Now, the additive Schwarz operator $T^{type}:V_h\rightarrow V_h$ is $$T^{type}=T^{(type)}_{0}+\sum^{N}_{i=1}T^{(i)},~~ type\in\{\RN{1},\RN{2}\}$$
and the problem \eqref{MainLinearForm} can be written as follow.
\begin{equation*}
T^{type}u_h^*=g^{type},~~type\in\{\RN{1},\RN{2}\},
\end{equation*}
where $g^{type}=g^{(type)}_{0}+\sum^{N}_{i=1}g^{(i)}$ with $g^{(type)}_{0}=T^{(type)}_{0}u_h^{*}$ and $g^{(i)}=T^{(i)}u_h^{*}$.

\section{Enrichment of the coarse space for the mortar discretization}\label{EnrichCorse}
In this section, we design two different coarse spaces for the additive average Schwarz method. To this end, we first denote the sets of nodal nodes of all mortar and nonmortar sides, and also all interior nodes of all subdomains $\Omega_i$, $i=1,\ldots,N$ by $\mathcal{N}_m$, $\mathcal{N}_s$, and $\mathcal{N}_i$, respectively. Now, the average interpolation operator  $I_0:V_h\rightarrow V_h$ for the mortar discretization as in \cite{Bjorstad-2003} has the following structure.
\begin{equation*}\label{interp}
I_0u^h(x)=
\begin{cases}
u^h_i(x)~~x\in     \mathcal{N}_m \cup \mathcal{N}_s ,\\
\overline{u}^h_i\hspace{0.7cm} x \in \mathcal{N}_i,
\end{cases}
i=1,\ldots,N,
\end{equation*}
where $\overline{u}^h_i$ is the average value of ${u}^h_i$ over $\Omega_i$, i.e.,
\begin{equation}\label{avg}
\overline{u}^h_i=\frac{1}{\mu_{(\delta,\gamma)}^i}\left(\sum_{\gamma_{m(i)}\subset\partial\Omega_{i}}\overline{u}^{h}_{\gamma_{m(i)}}~~+\sum_{\delta_{m(i)}\subset\partial\Omega_{i},\\ \gamma_{m(j)}=\delta_{m(i)}}\overline{u}^{h}_{\gamma_{m(j)}}\right),
\end{equation}
where
$$\overline{u}^h_{\gamma_{m(i)}}=\frac{1}{|\gamma_{m(i)}|}\int_{\gamma_{m(i)}}u^h_i ~ds$$ 
and $\mu_{(\delta,\gamma)}^i$ is the number of all mortar and nonmortar sides of $\Omega_i$, and $|\gamma_{m(i)}|$ is the length of $\gamma_{m(i)}$.
To express the interpolation operator $I_0$ in terms of the matrix form denoted by $\mathbf{R}_0$, let $\mathcal{N}_c$ be a set of all nodal nodes at the end of all nonmortar sides of $\Omega$ and let also  $\mathbf{I}_{(c)}\in \mathbb{R}^{N_c \times N_c}$ and  $\mathbf{I}_{(m)}\in \mathbb{R}^{N_m\times N_m}$ be identity matrices, where $N_c=\dim(\mathcal{N}_c)$ and $N_m=\dim(\mathcal{N}_m)$. Further, consider $\mathbf{H}=\text{diag}(\mathbf{H}_1,\ldots,\mathbf{H}_N) \in \mathbb{R}^{(N_c +N_m) \times N_i} $, where  $\mathbf{H}_i=(\overline{u}_i^h)_{l=1;~~j=1}^{4+n_m ~ n_i}$,  $N_i=\dim{(\mathcal{N}_i})$, and $n_i=\dim{(\imath_i})$. Hence
 \begin{equation*}
\mathbf{R}_{0} =\left[
  \begin{array}{c|c}
    \multicolumn{1}{c|}{\multirow{2}{*}{$\text{diag}(\mathbf{I}_{(c)}, \mathbf{I}_{(m)})$}} &  \multicolumn{1}{c}{\multirow{2}{*}{$\mathbf{H}$}}\\
    
    & 
  \end{array}
\right]   \in \mathbb{R}^{(N_{c} +N_{m}) \times N_{\Omega} }, ~ N_{\Omega} = N_c+N_m+N_i.
\end{equation*}
\begin{figure}[!ht]
\centering
\begin{overpic}[width=4.3cm,height=3.8cm]{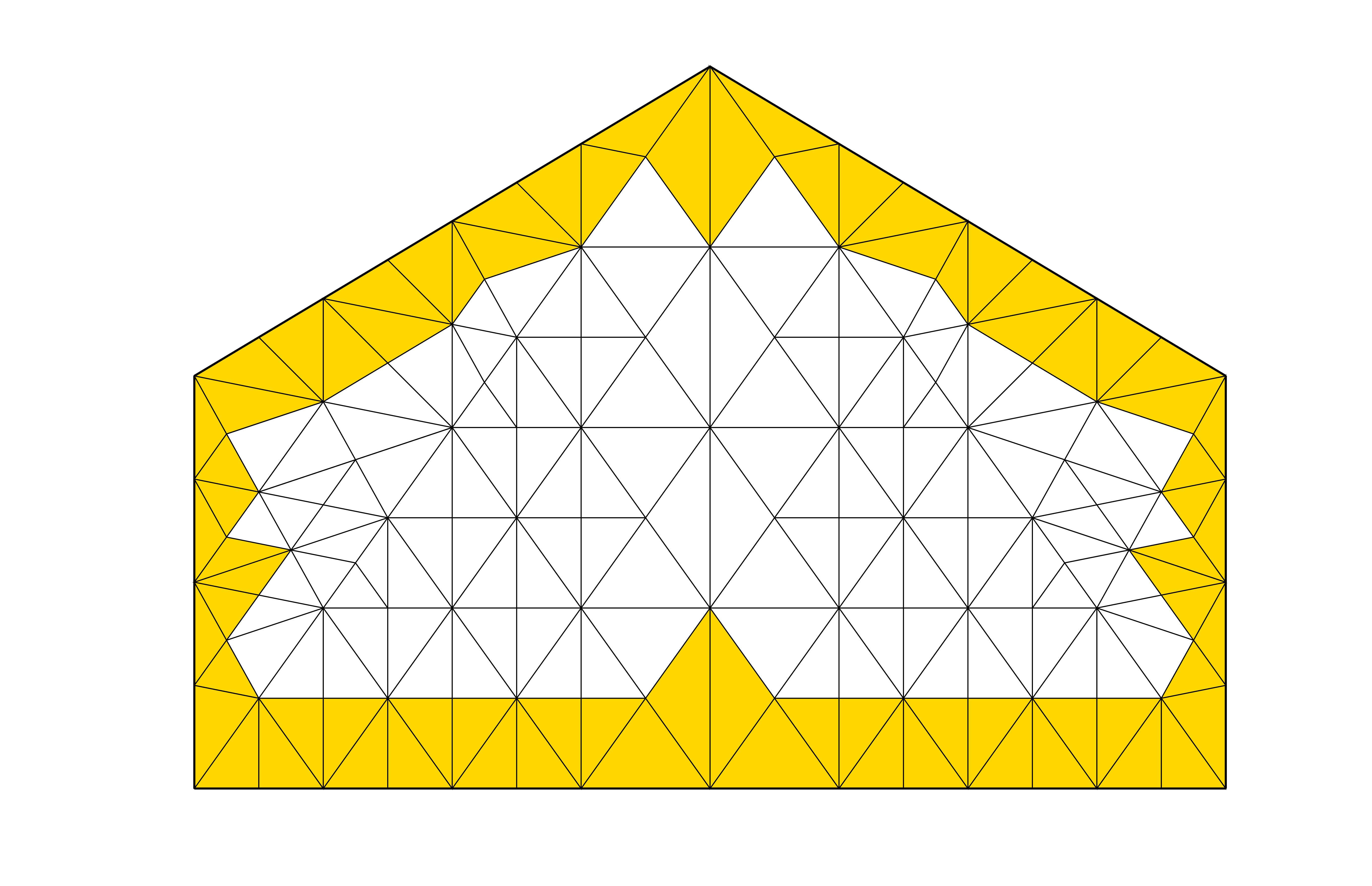}
 \put (60,22) {$\Omega_i$}
 \put (45,12) {$\Omega_i^{B}$}
 \end{overpic}
 \vspace{-0.3cm}
 \caption{The boundary layer $\Omega_{i}^{B} \subset \Omega_{i}$ is highlighted by the colorful triangles.}
 \label{MortarSFig1}
\end{figure}

To define two different coarse spaces, we first introduce the boundary layers $\Omega^B_i, ~i=1,\ldots,N$, where  $\Omega^B_i\subset\Omega_i$ is the sum of all triangles such as $\tau \in\mathcal{T}_i$ \vspace{-2cm}such that $~ \partial  \tau ~ \cap ~\partial\Omega_i ~ \neq ~ \phi ~$ (cf. Figure~\ref{MortarSFig1}).  We  then   define the following local\newpage \hspace{-0.5cm}minimums of the coefficients over the subdomains $\Omega_i$:
\begin{equation*}\label{minmaxcoefficients}
\underline{\alpha}_{i}:=\min_{x\in\Omega_i}\alpha(x),\qquad 
\underline{\alpha}_{i,B}:=\min_{x\in\Omega^B_i}\alpha(x),\quad 
i=1,\ldots,N. 
\end{equation*}

We now proceed to construct the second part of $V_{0}^{(type)}$, i.e.,  $W^{type}_{i},type\in\{\RN{1},\RN{2}\}, i=1,\ldots,N$ which are the spaces of adaptively chosen eigenfunctions of specially constructed generalized eigenvalue problems defined locally in each subdomain and extended by zero to the rest of the domain. Hence, the generalized eigenvalue problems is to find all eigen pairs $(\lambda^{i,type}_{k},\psi^{i,type}_{k})\in (\mathbb{R}_+,V^{(i)})$ such that
\begin{align}\label{GenEigPro}
\mathbf{A}_{\Omega_i}~\textbf{x} &= \lambda^{type} ~ (\mathbf{R}_i \mathbf{B}^{type} \mathbf{R}_i^T)~\textbf{x}, \quad i=1,\ldots,N,
\end{align}
where $\mathbf{B}^{type}(\cdot,\cdot)$, $type \in \{\RN{1},\RN{2}\}$ are the matrix representations of the following bilinear forms. 
\begin{align*}
b^{\RN{1}}(u,v)&:= \sum_{i=1}^N \int_{\Omega_i}\underline{\alpha}_i \nabla u \cdot \nabla v~dx, \quad u,v \in V_h,
\notag\\
b^{\RN{2}}(u,v)&:= \sum_{i=1}^{N} \left( \int_{\Omega_i^{B}}\underline{\alpha}_{i,B} \nabla u \cdot \nabla v~dx + \int_{\Omega_i \setminus
 \Omega_{i}^{B}} \alpha \nabla u \cdot \nabla v~dx\right), \quad u,v \in V_h.
\end{align*}

We also denote all eigenvalues of \eqref{GenEigPro} by
\begin{align*}
\mathbf{D}^{type} &= \text{diag}(\mathbf{D}_1^{type}, \mathbf{D}_2^{type}, \ldots, \mathbf{D}_N^{type}) \in \mathbb{R}^{ N_i \times  N_i}, ~ type \in \{ \RN{1}, \RN{2} \},
\end{align*}
where 
$$\mathbf{D}_i^{type} = \text{diag}( \lambda_1^{i,type}, \lambda_2^{i,type}, \ldots,  \lambda_{n_i}^{i,type})\in \mathbb{R}^{n_{i} \times n_{i}}, \quad i=1,\ldots,N,$$ such that
\begin{equation*}
   \lambda_1^{i,type}  \geq \lambda_1^{i,type} \geq \ldots \geq \lambda_{n_i}^{i,type} > 0.
\end{equation*}

We define
$$W^{type}_i:=\text{span}\{\psi^{i,type}_k\}^{n_i}_{k=1},~type\in\{\RN{1},\RN{2}\},i=1,\ldots,N$$
as the space of the eigenfunctions associated with the  eigenvalues $\lambda_{j}^{i,type}$. We also  correspond to these spectral spaces the  matrix forms denoted by $\mathbf{W}_i^{type},i=1,\ldots,N $ such that
$$
\mathbf{W}^{type} = \text{diag} (\mathbf{W}_1^{type}, \mathbf{W}_2^{type}, \ldots, \mathbf{W}_{N}^{type}), ~ type \in \{ \RN{1}, \RN{2} \}.
$$

Further, we use the notation $\mathbf{R}_{0}^{type}$ as the matrix representation of the operator $T^{(type)}_{0}$ which has the following structure.
\begin{equation*}
\mathbf{R}_{0}^{type}= 
\left[
  \begin{array}{c|c}
    \multicolumn{1}{c|}{\multirow{2}{*}{$\mathbf{R}_{0}$}} &  \multicolumn{1}{c}{\multirow{2}{*}{$O$}}\\
    
     \\
    \hline
    \multicolumn{1}{c|}{\multirow{2}{*}{ $O$}}
    & \multicolumn{1}{c}{\multirow{2}{*}{  $ \mathbf{W}^{type}$ }}\\
    & \\
  \end{array}
\right],\quad \mathbf{R}_{0}^{type} \in \mathbb{R}^{N_{\Omega} \times (N_{\Omega}+N_i)}, ~ type \in \{ \RN{1}, \RN{2} \}.
\end{equation*}

\subsection{The additive average Schwarz preconditioner }

The main aim of this section is to introduce new  enrichment preconditioners denoted by $\mathbf{B}_{E}^{type}, type \in \{ \RN{1}, \RN{2} \}$. The idea is based on expressing these preconditioners in terms of a combination of 
$\mathbf{B}_C^{type}$, $type \in \{  \RN{1},\RN{2} \}$ and $\sum_{i=1}^{N} \mathbf{R}_i^{T} \mathbf{A}_{\Omega_i}^{-1}\mathbf{R}_i,$ where $\mathbf{B}_C^{type} = (\mathbf{R}_{0}^{type}) ^{T} \Big( \mathbf{R}_{0}^{type} \mathbf{A}_{N} (\mathbf{R}_{0}^{type}) ^{T} \Big)^{-1} \mathbf{R}_{0}^{type}$ and
\begin{equation*}
\mathbf{A}_{N}=
\left[
  \begin{array}{c|c}
    \multicolumn{1}{c|}{\multirow{2}{*}{$\mathbf{A}^{(11)}_{N}$}} &  \multicolumn{1}{c}{\multirow{2}{*}{$\mathbf{A}^{(12)}_{N}$}}\\
     \\
    \hline
    \multicolumn{1}{c|}{\multirow{2}{*}{ $\mathbf{A}^{(21)}_{N}$}}
    & \multicolumn{1}{c}{\multirow{2}{*}{ $\mathbf{A}^{(22)}_{N}$ }}\\
    & \\
  \end{array}
\right], \quad \mathbf{A}_{N} \in \mathbb{R}^{(N_{\Omega}+N_i) \times (N_{\Omega}+N_i)}
\end{equation*}
such that
\begin{align*}
& \begin{minipage}{33em} $\mathbf{A}_{N}^{(11)}:=$ square  submatrix of $\mathbf{A}$ whose rows and columns are corresponding to the nodes belong to $\mathcal{N}_c \bigcup \mathcal{N}_m \bigcup \mathcal{N}_i$, \end{minipage}\\
    & \begin{minipage}{33em} $\mathbf{A}_{N}^{(12)}  :=$ rectangular submatrix of $\mathbf{A}$ whose rows are corresponding to all  nodes belong to $\mathcal{N}_c \bigcup \mathcal{N}_m \bigcup \mathcal{N}_i $, meanwhile whose columns are corresponding to the nodes belong to $\mathcal{N}_i$,\\ \end{minipage}\\
                              &  \begin{minipage}{33em} $\mathbf{A}_{N}^{(21)}= \left( \mathbf{A}_{N}^{(12)} \right)^{T}$,\end{minipage} \\
    & \begin{minipage}{33em} $\mathbf{A}_{N}^{(22)} :=$ square submatrix of $\mathbf{A}$ whose rows and columns are corresponding to the nodes belong to $\mathcal{N}_i$. Indeed, $\mathbf{A}_{N}^{(22)}= \text{diag}(\mathbf{A}_{\Omega_1}, \mathbf{A}_{\Omega_2}, \ldots, \mathbf{A}_{\Omega_N} )$.\end{minipage}
\end{align*}

Hence
\begin{equation*}
\mathbf{R}_{0}^{type} \mathbf{A}_{N} (\mathbf{R}_{0}^{type}) ^{T}=
\left[
  \begin{array}{c|c}
    \multicolumn{1}{c|}{\multirow{2}{*}{$\mathbf{R}_{0}\mathbf{A}^{(11)}_{N}\mathbf{R}_{0}^{T}$}} &  \multicolumn{1}{c}{\multirow{2}{*}{$\mathbf{R}_{0}\mathbf{A}^{(12)}_{N}(\mathbf{W}^{type})^{T}$}}\\
    
     \\
    \hline
    \multicolumn{1}{c|}{\multirow{2}{*}{ $\Big(\mathbf{R}_{0}\mathbf{A}^{(12)}_{N}(\mathbf{W}^{type})^{T} \Big)^{T}$}}
    & \multicolumn{1}{c}{\multirow{2}{*}{$   \mathbf{W}^{type}\mathbf{A}^{(22)}_{N}(\mathbf{W}^{type})^{T} $ }}\\
    & \\
  \end{array}
\right].
\end{equation*}

Since $\mathbf{W}^{type}\mathbf{A}^{(22)}_{N}(\mathbf{W}^{type})^{T}=\mathbf{D}^{type}, type \in\{ \RN{1}, \RN{2} \}$ and it is the nonsingular matrix, we can easily compute $ \mathbf{B}_C^{type}$. To this end and for the sake of simplicity, we first assume
\begin{align*}
\mathbf{G}^{type}& = \mathbf{R}_{0}\mathbf{A}^{(12)}_{N}(\mathbf{W}^{type})^{T},\\
\mathbf{S}^{type}&=\mathbf{R}_{0}\mathbf{A}^{(11)}_{N}\mathbf{R}_{0}^{T} - \mathbf{G}^{type} \big( \mathbf{D}^{type} \big)^{-1} \big(\mathbf{G}^{type} \big)^{T}.
\end{align*}

Hence, the inverse of $2 \times 2$ block matrix  $\mathbf{R}_{0}^{type} \mathbf{A}_{N} (\mathbf{R}_{0}^{type}) ^{T}$ implies
\begin{align*}
 \mathbf{B}^{type}_{C(11)} & = \mathbf{R}_{0}^{T}(\mathbf{S}^{type})^{-1} \mathbf{R}_{0},\\
\mathbf{B}^{type}_{C(12)} & = - \mathbf{R}_{0}^{T}  (\mathbf{S}^{type})^{-1}  \mathbf{G}^{type} (\mathbf{D}^{type})^{-1}\mathbf{W}^{type},\\
\mathbf{B}^{type}_{C(21)} & = - \Big( \mathbf{B}^{type}_{C(12)} \Big)^T,\\
\mathbf{B}^{type}_{C(22)} & = (\mathbf{W}^{type})^{T} \Big(   (\mathbf{D}^{type})^{-1} +  \big(  \mathbf{G}^{type} (\mathbf{D}^{type})^{-1}\big)^{T} (\mathbf{S}^{type})^{-1} \mathbf{G}^{type} (\mathbf{D}^{type})^{-1}  \Big) \mathbf{W}^{type}.
\end{align*}

Due to different sizes of the submatrices of $\mathbf{B}_{C}^{type}$, we can use the restriction matrices  $\mathbf{R}_i, i=1,\ldots,N$ to obtain new matrices not only with identical sizes, but also with their  previous performances. Let $\mathbf{R}_{r} = \text{diag}(\mathbf{R}_1, \mathbf{R}_2, \ldots, \mathbf{R}_N)$, $ \mathbf{B}_{0}^{type} = \mathbf{B}^{type}_{C(11)} + \mathbf{R}_{r}^{T} \mathbf{B}^{type}_{C(21)}$ and $\mathbf{B}_{00}^{type} = \mathbf{B}^{type}_{C(12)} \mathbf{R}_{r} +\mathbf{R}_{r}^{T}  \mathbf{B}^{type}_{C(22)} \mathbf{R}_{r}.
$ Hence, $\mathbf{B}_{E}^{type}, type \in \{ \RN{1},\RN{2} \}$  as the enrichment preconditioners can be written explicitly in the following matrix form. 
\begin{align*}
\mathbf{B}_{E}^{type} & = \mathbf{B}^{type}_{0} + \mathbf{B}_{00}^{type}+\sum_{i=1}^{N} \mathbf{R}_i^{T} \mathbf{A}_{\Omega_i}^{-1}\mathbf{R}_i\\
& =\mathbf{B}^{type}_{0} + \left( \mathbf{R}_{r}^{T} -  \mathbf{B}^{type}_{0} \mathbf{A}_{N}^{(12)}  \right)  (\mathbf{W}^{type})^{T} (\mathbf{D}^{type})^{-1} \mathbf{W}^{type} \mathbf{R}_{r} + \sum_{i=1}^{N} \mathbf{R}_i^{T} \mathbf{A}_{\Omega_i}^{-1}\mathbf{R}_i.
\end{align*}

To estimate an upper bound for the condition number of   $\mathbf{B}_{E}^{type}$ in the next section, we define the $b^{type}_i(\cdot,\cdot)$-orthogonal projection operator $\pi^{type}_i:V^{(i)}\rightarrow V^{(i)}$ as
$$\pi^{type}_iv^h=\sum^{N_i}_{k=1}b^{type}_i(v^h,\psi^{i,type}_k)\psi^{i,type}_k,\quad v^h \in V_h, \quad type\in\{\RN{1},\RN{2} \}.$$ 

For any $u^h\in V_h$, we consider the function $w^h=u^h-I_0u^h\in V_h$ with zero value both on $\partial\Omega_i$ and on the rest of the domain $\Omega$. In addition, we define  $I^{type}_0:V_h\rightarrow V_0^{(type)}$ as
\begin{equation*}
I^{type}_0 u^h=I_0 u^h+\sum^{N}_{i=1}\pi^{type}_i w^h,\qquad u^h \in V_h, \quad type\in \{\RN{1},\RN{2}\}.
\end{equation*}

\section{On the estimation of the condition number bound}\label{CondNumb}

This section is devoted to obtaining condition number estimates of the additive average Schwarz preconditioners, which are defined in the previous section for $type \in \{\RN{1},\RN{2}\}$.  The proof is based on the standard additive Schwarz framework, where three assumptions to be held, see \cite[p.~155]{Smith-book-1996}. Here, we need only to show that there exists the stable splitting for all $u \in V_{h}$ as in Theorem \ref{SplittingTheorem}. Note that the other assumptions hold since there is no overlapping, and we also use the exact bilinear forms. For simplicity, throughout this section, we use $C$ as a positive constant, which is independent of the mesh sizes. We add a notation or an index to $C$ if we want to emphasize a special constant. Furthermore, we use the notation $\lesssim$ to remove any constants except the mesh sizes in the inequalities.

Let $M^{type}_i$ be a given number such that $0\leq M^{type}_i < n_i$ and $\lambda^{i,type}_{M^{type}_i+1}<\lambda^{i,type}_{M^{type}_i}$. We define
\begin{equation*}
\widetilde{W}^{type}_i:=\text{span}\{\psi^{i,type}_k\}^{M_i^{type}}_{k=1},~type\in\{\RN{1},\RN{2}\}
\end{equation*}
and 
\begin{equation*}
\widetilde{V}_0^{(type)}=I_0V_h+\sum\nolimits^{N}_{i=1}\widetilde{W}^{type}_{i},\quad type\in\{\RN{1},\RN{2}\}.\end{equation*}

For the analysis of the  additive average Schwarz method, we define $\widetilde{I}_0^{type}:V_h \rightarrow \widetilde{V}_0^{(type)}$ as 
\begin{equation*}
\widetilde{I}^{type}_0 u^h=I_0 u^h+\sum\nolimits^{N}_{i=1}\widetilde{\pi}^{type}_i w^h,\quad u^h \in V_h, \quad type\in \{\RN{1},\RN{2}\}, 
\end{equation*}
where 
\begin{equation*}
\widetilde{\pi}^{type}_iv^h=\sum\nolimits^{M_i^{type}}_{k=1}b^{type}_i(v^h,\psi^{i,type}_k)\psi^{i,type}_k, \quad type\in\{\RN{1},\RN{2}\}.
\end{equation*}

To prove our main result, we need the following lemma. 
\begin{lem}\label{lemmaOnL2}
For all $u \in V^{(i)}, i=1,\ldots,N$,
\begin{align*}
    |u - \widetilde{\pi}_{i}^{\RN{1}} u |^2_{H^1(\Omega_i),\alpha} & \leq C \lambda_{M_i^{\RN{1}}+1}^{i,\RN{1}} \| \underline{\alpha}_i^{1/2} \nabla u \|^2_{L^2(\Omega_i)},\\
    |u - \widetilde{\pi}_{i}^{\RN{2}} u |^2_{H^1(\Omega_i),\alpha} & \leq C \lambda_{M_i^{\RN{2}}+1}^{i,\RN{2}} \Big( \| \underline{\alpha}_{iB}^{1/2} \nabla u\|^2_{L^2(\Omega_{i}^B)} + \| \alpha_{i}^{1/2} \nabla u\|^2_{L^2(\Omega_i \setminus \Omega_{i}^B)}    \Big),
\end{align*}
 where ${|\cdot|}^2_{H^1(\Omega_i),\alpha}=a_i(\cdot,\cdot)$.
\begin{proof}
\rm 
We first express any  $u \in V^{(i)}$ uniquely in terms of the  eigenfunctions, i.e., $u = \sum_{k=1}^{n_i} b_i^{type}(u,\psi_{k}^{i,type}) \psi_{k}^{i,type}$.  Hence 
\begin{equation*}
    u- \widetilde{\pi}_i^{type}u = \sum\nolimits_{k=M_{i}^{type}+1}^{n_i} b_i^{type}(u,\psi_{k}^{i,type}) \psi_{k}^{i,type}, \quad type \in \{\RN{1},\RN{2}\}
\end{equation*}
and we have
\begin{equation*}
    | u- \widetilde{\pi}_i^{type}u |_{H^1(\Omega_i),\alpha}^2 \leq   \sum\nolimits_{k=M_{i}^{type}+1}^{n_i}\left( b_i^{type}(u,\psi_{k}^{i,type} ) \right)^2 | \psi_{k}^{i,type} | _{H^1(\Omega_i),\alpha}^2.
\end{equation*}

Furthermore, $(\mathbf{W}_i^{type})^T \mathbf{A}_{\Omega_i}\mathbf{W}_i^{type} = \mathbf{D}_i^{type}$ is equivalent to 
\begin{align*}
    |\nabla \psi_{k}^{i,type}|^2_{H^1(\Omega_i),\alpha} & = \lambda_{k}^{i,type}, \quad k=1,\ldots, n_i \quad i=1,\ldots,N, \quad type \in \{\RN{1},\RN{2}\}.
\end{align*}

Therefore 
\begin{equation}\label{GenEiq}
    | u- \widetilde{\pi}_i^{type}u |_{H^1(\Omega_i),\alpha}^2 \leq  \sum\nolimits_{k=M_{i}^{type}+1}^{n_i}\left( b_i^{type}(u,\psi_{k}^{i,type} ) \right)^2\lambda_{k}^{i,type}, \quad type \in \{\RN{1},\RN{2}\}.
\end{equation}

Using the Schwarz inequality for $type=\RN{1}$, yields

\begin{equation*}
     | u- \widetilde{\pi}_i^{\RN{1}}u |_{H^1(\Omega_i),\alpha}^2 \leq \|\underline{\alpha}_i^{1/2}\nabla u \|^2_{L^2(\Omega_i)}  \sum\nolimits_{k=M_{i}^{\RN{1}}+1}^{n_i}  \|\underline{\alpha}_i^{1/2}\nabla \psi_{k}^{i,\RN{1}} \|^2_{L^2(\Omega_i)}\lambda_{k}^{i,\RN{1}}. 
\end{equation*}

Since $(\mathbf{W}_i^{type})^T (\mathbf{R}_i \mathbf{B}^{type} \mathbf{R}_i^T) \mathbf{W}_i^{type} = \mathbf{I}, ~type \in \{ \RN{1}, \RN{2} \}$ and $i=1,\ldots,N$ we have equivalently 
\begin{equation}\label{PropertyEig}\|\underline{\alpha}_i^{1/2}\nabla \psi_{k}^{i,\RN{1}} \|^2_{L^2(\Omega_i)}=1, \quad  \|\underline{\alpha}_{iB}^{1/2} \nabla \psi_{k}^{i,\RN{2}} \|^2_{L^2(\Omega_{i}^B)} + \|\alpha_{i}^{1/2} \nabla \psi_{k}^{i,\RN{2}} \|^2_{L^2(\Omega_i \setminus \Omega_{i}^B)} =1.
\end{equation}

Then 
\begin{equation*}
    | u- \widetilde{\pi}_i^{\RN{1}}u |_{H^1(\Omega_i),\alpha}^2 \leq C  \lambda_{M_i^{\RN{1}}+1}^{i,\RN{1}} \|\underline{\alpha}_i^{1/2} \nabla u \|^2_{L^2(\Omega_i)}. 
\end{equation*}

For $type=\RN{2}$, from \eqref{GenEiq} we get
\begin{align*}
    | u- \widetilde{\pi}_i^{\RN{2}}u |_{H^1(\Omega_i),\alpha}^2 & \leq C  \sum\nolimits_{k=M_{i}^{\RN{2}}+1}^{n_i}  \Big(   \| \underline{\alpha}_{iB}^{1/2} \nabla u\|^2_{L^2(\Omega_{i}^B)} \|\underline{\alpha}_{iB}^{1/2} \nabla \psi_{k}^{i,\RN{2}} \|^2_{L^2(\Omega_{i}^B)} \notag\\ 
    & \hspace{0.8cm} +  \| \alpha_{i}^{1/2} \nabla u\|^2_{L^2(\Omega_i \setminus \Omega_{i}^B)} \|\alpha_{i}^{1/2} \nabla \psi_{k}^{i,\RN{2}} \|^2_{L^2(\Omega_i \setminus \Omega_{i}^B)} \Big) \lambda_{k}^{i,\RN{2}}.
\end{align*}

To complete the proof, it is sufficient to use \eqref{PropertyEig}.  
\end{proof}
\end{lem}

\begin{thm}\label{SplittingTheorem}
For all $u^h\in V_h$ the following results hold:
\begin{equation*}
a(\widetilde{I}^{type}_0 u^h , \widetilde{I}^{type}_0 u^h) \lesssim \max_{i} \lambda_{M_i^{type}+1}^{i,type} ~\frac{H}{h} a(u^h,u^h), \quad type \in\{\RN{1},\RN{2}\}, \quad i=1,\ldots,N.
\end{equation*}
\begin{proof}
\rm 
We first consider the following splitting.
\begin{equation*}
    u^h = \widetilde{I}^{type}_0 u^h +\sum_{i=1}^{N}u^i,
\end{equation*}
where $u^i=(0,\ldots,0,w_i,0,\ldots,0) \in V^{(i)}$ and $w_i=(u^h - \widetilde{I}^{type}_0 u^h)|_{\Omega_i}$. Consequently, we get
\begin{align*}
a(\widetilde{I}^{type}_0 u^h , \widetilde{I}^{type}_0 u^h ) &\preceq a(u^h,u^h) + a(u^h - \widetilde{I}^{type}_0 u^h , u^h - \widetilde{I}^{type}_0 u^h),
\end{align*}
where
\begin{align*} 
    a(u^h - \widetilde{I}^{type}_0 u^h ,u^h - \widetilde{I}^{type}_0 u^h ) &= \sum_{i=1}^{N}a(u^i,u^i)
    =\sum_{i=1}^{N}a_i(u^h - \widetilde{I}^{type}_0 u^h ,u^h - \widetilde{I}^{type}_0 u^h )\\
    & = \sum_{i=1}^{N} a_i (    w - \widetilde{\pi}_i^{type}w , w - \widetilde{\pi}_i^{type}w  ), \quad type \in \{ \RN{1}, \RN{2}  \}.
\end{align*}

We first use Lemma \ref{lemmaOnL2} for the $type= \RN{1}$. Then, we have
\begin{align*}
     a_i (    w - \widetilde{\pi}_i^{\RN{1}}w , w - \widetilde{\pi}_i^{\RN{1}}w  ) &\leq C \lambda_{M_i^{\RN{1}}+1}^{i,\RN{1}} \underline{\alpha}_i{||\nabla(u^h-I_0u^h)||}^{2}_{L^{2}(\Omega_i)}\\
& \leq  C \lambda_{M_i^{\RN{1}}+1}^{i,\RN{1}} \underline{\alpha}_i \left( \frac{H}{h} ||\nabla u^h||^{2}_{L^{2} (\Omega_i)}  + ||\nabla I_0 u^h||^{2}_{L^{2}(\Omega_i)}  \right).
\end{align*}

Further, from \cite[pp.~8-9]{Bjorstad-2003} we get similarly 
\begin{equation}\label{FresultI}
    \underline{\alpha}_i ||\nabla I_0 u^h||^{2}_{L^{2}(\Omega_i)} \leq C \frac{H}{h}a_i(u^h,u^h).
\end{equation}

For $type= \RN{1}$, the proof is completed by summing \eqref{FresultI} over $i=1,\ldots,N$. We now proceed to prove a similar result for $type= \RN{2}$. Due to definition of $I_0 u_i^h$, we use this fact that $\nabla I_0 u_i^h$ is equal to zero on each triangle $\tau \not\in \Omega_i^{B}$, $i=1,\ldots,N$. Hence, we use Lemma \ref{lemmaOnL2} to get the first inequality as follow.

\begin{align*}
    a_i (    w & - \widetilde{\pi}_i^{\RN{2}}w , w - \widetilde{\pi}_i^{\RN{2}}w  )\notag\\ & \leq C \lambda_{M_i^{\RN{2}}+1}^{i,\RN{2}}
    \Big( {||\alpha_i^{1/2}\nabla(u^h-I_0u^h)||}^{2}_{L^{2}(\Omega_i\setminus\Omega_i^B)}+\underline{\alpha}_{i,B}{||\nabla(u^h-I_0u^h)||}^{2}_{L^{2}(\Omega_i^B)} \Big) \\
& \leq C \lambda_{M_i^{\RN{2}}+1}^{i,\RN{2}}
    \Big( {||\alpha^{1/2}_{i}\nabla u^h||}^{2}_{L^{2}(\Omega_i\setminus\Omega^B_{i})}+\underline{\alpha}_{i,B}
{||\nabla(u^h-I_0u^h)||}^{2}_{L^{2}(\Omega_i^B)} \Big) \\
&\leq C \lambda_{M_i^{\RN{2}}+1}^{i,\RN{2}}
    \bigg( {||\alpha^{1/2}_i\nabla u^h||}^{2}_{L^{2}(\Omega_i)}+\underline{\alpha}_{i,B} \left( ||\nabla u^h ||_{L^{2}(\Omega_i^B)}^{2} +  || \nabla I_0u^h||_{L^{2}(\Omega_i^B)}^2 \right) \bigg)\\
&\leq C \lambda_{M_i^{\RN{2}}+1}^{i,\RN{2}}
    \left( a_i(u^h,u^h) + {||\underline{\alpha}_{i}^{1/2}\nabla  u^h||}^{2}_{L^{2}(\Omega_i)} +  \underline{\alpha}_{i,B}{||\nabla  I_0u^h||}^{2}_{L^{2}(\Omega_i^B)} \right) \notag \\
& \leq C \lambda_{M_i^{\RN{2}}+1}^{i,\RN{2}}
    \left( a_i(u^h,u^h) + \underline{\alpha}_{i,B}{||\nabla  I_0u^h||}^{2}_{L^{2}(\Omega_i^B)}\right).
\end{align*}

Using the inverse inequality and the definition of the operator $I_0$, implies 

\begin{align}
{||\nabla I_0u^h||}^{2}_{L^{2}(\Omega_i^B)}
&=\sum_{\tau \in \overline{\Omega}^B_i}{||\nabla I_0u^h||}^{2}_{L^{2}(\tau)}\notag\\
&=\sum_{\tau \in \overline{\Omega}^B_i}{||\nabla(I_0u^h -\overline{u}^h)||}^{2}_{L^{2}(\tau)}\notag\\
&\leq C\sum_{\tau \in \overline{\Omega}^B_i}h^{-2}_{\tau}{||I_0u^h -\overline{u}^h||}^{2}_{L^{2}(\tau)}\notag\\
&\leq C\sum_{x\in\partial\Omega_{ih}}(u^h_i(x)-\overline{u}^h_i)^2.\label{ApInter}
\end{align}

In what follow, we need to use the following Poincar\'e inequality \cite{Toselli-2005}.  
\begin{equation}\label{PoincareInequality}
    ||g||_{L^2(\Omega_i)}^2 \leq C_1 |g|_{H^1(\Omega_i)}^2 + C_2 \left( \int_{\Omega_i} g ~ dx \right)^2, \quad g \in H^1(\Omega_i),
\end{equation}
where $C_1$ and $C_2$ are two positive constants independent of the mesh size of $\Omega$. Consider the function
$g=u_i^h - c$, where $c=\frac{1}{\mathrm{meas}(\Omega_i)}\int_{\Omega_i}u_i^h~ dx$. Now, (\ref{PoincareInequality})   implies
\begin{equation}\label{resultedfromPoinIneq}
||u_i^h-c||_{L^2(\Omega_i)}^2 \leq C_1 |u_i^h|_{H^1(\Omega_i)}^2.
\end{equation}

Furthermore, we use (\ref{avg}) and this property that all boundary elements of subdomains $\Omega_i,i=1,\ldots,N$ are quasi-uniform which implies that the number of all nodes belonging to $\partial \Omega_i$, denoted by $M_i$, is of the order $H_i/h_i$, to get 
\begin{align}
M_i\left(\overline{u^h_i-c}\right)^2 & \leq CM_i\frac{1}{\mu_{(\delta,\gamma)}^i} \sum_{\Gamma\subset\partial\Omega_i}
\frac{1}{\Gamma}{||u^h_i-c||}^{2}_{L^{2}(\Gamma)}\notag\\
&\leq Ch^{-1}_i{||u^h_i-c||}^{2}_{L^{2}(\partial\Omega_i)},\label{ApAvg}
\end{align}
where $\Gamma$ is $\gamma_{m(i)}$ and $\delta_{m(i)}$.

Now, we use $g$ and  (\ref{ApAvg}) to estimate an upper bound for  (\ref{ApInter}) as follows. 
\begin{align}
    {||\nabla I_0u^h||}^{2}_{L^{2}(\Omega_i^B)}
& \leq C \sum_{x\in\partial\Omega_{ih}}
\left(g(x)-\overline{g}\right)^2  \notag\\
&\leq Ch^{-1}_i{||g||}^{2}_{L^{2}(\partial\Omega_i)}\notag\\
& \leq C h_i^{-1}H_i ||\hat{g}||^{2}_{L^{2}(\partial \hat{\Omega})},\notag
\end{align}
where $\hat{\Omega}$ is the reference element of unit diameter. From this inequality we deduce 
\begin{align*}
||\nabla I_0u^h||^{2}_{L^{2}(\Omega_i^B)} & \leq Ch_i^{-1}H_i\left( |\hat{g}|_{H^1(\hat{\Omega})}^2 + || \hat{g} ||^2_{L^{2}(\hat{\Omega })} \right)\notag\\
& \leq C h_i^{-1}\left(H_i{|g|}^{2}_{H^{1}(\Omega_i)}+H^{-1}_i{||g||}^{2}_{L^{2}(\Omega_i)}\right)\notag\\
& \leq C \frac{H_i}{h_i}a_i(u^h,u^h),
\end{align*}
by the trace inequality, Theorem 3.1.2 in \cite{Ciarlet-1978} and \eqref{resultedfromPoinIneq}. To complete the proof, it suffices to take a summation over $i=1,\ldots,N$. 
\end{proof}
\end{thm}

To estimate upper bounds of the condition numbers of $\mathbf{B}_{E}^{type} \mathbf{A}$, $type \in \{\RN{1},\RN{2} \}$ it suffices to use Lemma 3 in \cite[pp.~156-158]{Smith-book-1996} and Theorem \ref{SplittingTheorem}. 

\begin{thm}
The condition numbers of the enriched additive average Schwarz preconditioners are bounded by
\begin{equation*}
    \kappa (\mathbf{B}_{E}^{type} \mathbf{A}) \leq C \left( \frac{H}{h} \right) \max_{i} \lambda_{M_i^{type}+1}^{i,type}, \quad type \in\{\RN{1},\RN{2}\}, \quad i=1,\ldots,N,
\end{equation*}
where $C$ is a positive constant. 
\end{thm}

\section{Numerical results}\label{SNumericalResults}
\begin{figure}[!ht]
\centering
\vspace{-.8cm}
\setlength{\unitlength}{0.1\textwidth}
\begin{picture}(15,6.5)
\put(1.3,-0.7){ \hspace{0cm} \includegraphics[width=8.cm]{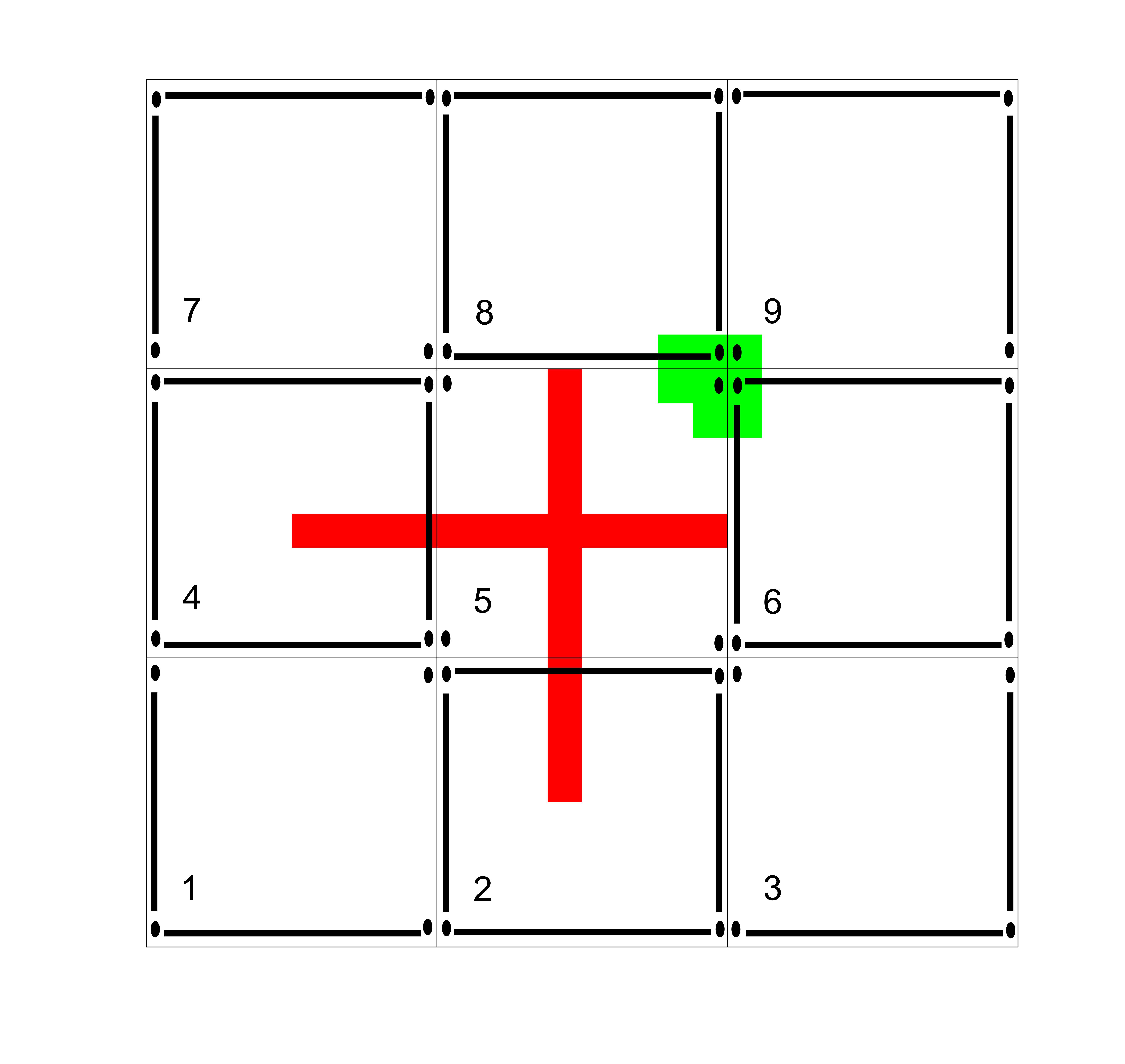}}
  \end{picture}
  \caption{The decomposition of $\Omega=[0,1]\times[0,1]$ into $3 \times 3$ subdomains (numbered from 1 to 9) and the locations of $\alpha(x)$  consisting of  $\alpha_b$, $\alpha_c$ and $\alpha_i$  in the white, red and green areas, respectively. Further, all mortar sides are denoted by the thick and black line segments.} 
  \label{first alpha0}
\end{figure}
\begin{figure}[!ht]
\centering
\setlength{\unitlength}{0.1\textwidth}
\vspace{-.8cm}
\begin{picture}(15,6.5)
\put(-1.05,-0.7){ \hspace{0cm} \includegraphics[width=13.45cm]{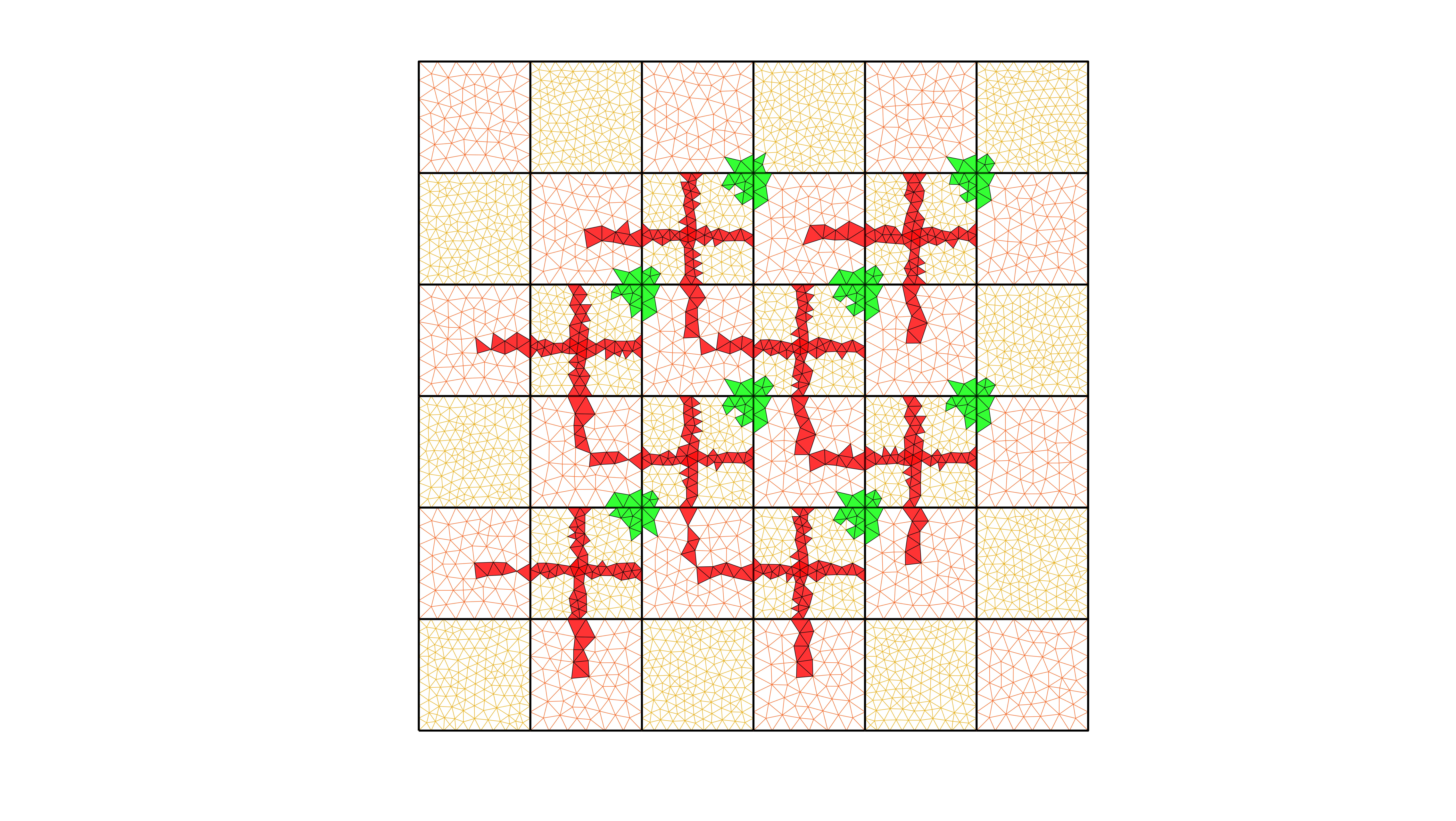}}
  \end{picture}
  \caption{The distribution of all jumps in $\alpha(x)$ following the  extended pattern in Figure~\ref{first alpha0}, where $\Omega=[0,1] \times [0,1]$ is divided into $6 \times 6$ subdomains. The nonmatching discretization parameters are  $H=1/6$, $h=1/36$ and $h^{*}=1/54$. Further, all mortar sides are the coarse mortar.}
  \label{first alpha}
\end{figure}

In this section, numerical results confirm the additive average Schwarz method's validity and efficiency with adaptive enrichment, where the jumps in the coefficient $\alpha(x)$ in \eqref{eq:general_model} are very large and even change rapidly. Those jumps might be occurred inside of the subdomains, or even across the subdomain boundaries. To have a complicated distribution of jumps, we use the following pattern \cite{Marcinkowski-2018}, depicted in Figure~\ref{first alpha0}, which is a periodic pattern when the number of subdomains are increased (cf. Figure~\ref{first alpha} ). We also consider background channels, crossing channels and corner channels denoted by $\alpha_b$, $\alpha_i$ and  $\alpha_c$, respectively. For different values of $\alpha_b$, $\alpha_c$ and $\alpha_i$, our test problem has the right-hand side function $f(x,y) = 2 \pi^{2} \sin(\pi x)\sin(\pi y)$ defined in $\Omega=[0,1]\times[0,1]$.

All presented numerical results are based on the nonmatching triangulations across the subdomains interfaces, where $H\in \{1/6,1/9\}$, $h\in \{    1/36, 1/54\}$, and $h^{*}\in\{ 1/54,1/81\}$ (cf. Figure~\ref{first alpha}). We use the similar locations for $\alpha_c$ and $\alpha_i$ as the periodic patterns and the similar mortar sides, as depicted in  Figure~\ref{first alpha0}, 
for $N \geq 9$. For instance, see Figure~\ref{first alpha}, where $N=36$. Further, in what follows, we use two terms, i.e., coarse mortar and fine mortar to distinguish between the coarse and fine triangulations connected to the mortar sides. To withdraw numerical results, we use the additive average Schwarz method to produce the enrichment  preconditioners. We estimate the condition number of such preconditioners for $type \in \{ \RN{1},\RN{2}\}$. Moreover, the iteration numbers in all tables come from the preconditioned conjugate gradient method, with the tolerance $5e-6$,  based on the produced preconditioners.

\begin{table}[!ht]
\begin{center}
\begin{tabular}{|c|cccc c c c cccccc | }
\cline{1-14}
\multicolumn{2}{|c|}{\multirow{3}{*}{ADD}} & \multicolumn{12}{c|}{ \tiny{ Different Values for} $\alpha_b$, $\alpha_c$ and $\alpha_i$ }\\
\cline{3-14}
\multicolumn{2}{|c|}{} & \multicolumn{6}{|c|}{ $1$, $1\text{e}3$,  $1\text{e}4$} & \multicolumn{6}{|c|}{$1$,  $1\text{e}4$,  $1\text{e}6$ } \\
\cline{3-14}
\multicolumn{2}{|c|}{}& \multicolumn{3}{|c|}{\tiny Coarse Mortar }& \multicolumn{3}{|c|}{ \tiny Fine Mortar } & \multicolumn{3}{|c|}{ \tiny Coarse Mortar}  &  \multicolumn{3}{|c|}{\tiny Fine Mortar }\\
\hline
\multirow{6}{*}{\begin{sideways} {\tiny	Number of Subdomains \text{ }} \end{sideways}}  & \multicolumn{1}{c|}{ \multirow{4}{*}{$6 \times 6$} } & \multicolumn{3}{c}{ \multirow{3}{*}{$6.07\text{e}1 $ }}&\multicolumn{3}{c|}{ \multirow{3}{*}{ $4.25\text{e}1$}}& \multicolumn{3}{c}{\multirow{3}{*}{$5.99\text{e}1$ }} &  \multicolumn{3}{c|}{\multirow{3}{*}{ $4.24\text{e}1$}} \\ 
&\multicolumn{1}{c|}{} & & & & & & \multicolumn{1}{c|}{}& & & & & &\\
& \multicolumn{1}{c|}{}& \multicolumn{3}{c}{\multirow{1}{*}{$(52)$}}    & \multicolumn{3}{c|}{\multirow{1}{*}{$(45)$}} & \multicolumn{3}{c}{\multirow{1}{*}{$(56)$}} & \multicolumn{3}{c|}{\multirow{1}{*}{$(49)$}}   \\
& \multicolumn{1}{c|}{\multirow{6}{*}{$9\times9$}} & \multicolumn{3}{c}{\multirow{5}{*}{$7.99\text{e}1$ }}& \multicolumn{3}{c|}{\multirow{5}{*}{$5.56\text{e}1$}} & \multicolumn{3}{c}{\multirow{5}{*}{$7.69\text{e}1$ }}  &  \multicolumn{3}{c|}{\multirow{5}{*}{$5.51\text{e}1$}}\\
\cline{2-14} 
& \multicolumn{1}{c|}{}& \multicolumn{3}{c}{\multirow{5}{*}{$(55)$}}    &   \multicolumn{3}{c|}{\multirow{5}{*}{$(50)$}}  &   \multicolumn{3}{c}{\multirow{5}{*}{$(58)$}}   &    \multicolumn{3}{c|}{\multirow{5}{*}{$(52)$}}\\
& \multicolumn{1}{c|}{}& & & &  & & \multicolumn{1}{c|}{}& & & & & & \\
& \multicolumn{1}{c|}{}& & & &  & & \multicolumn{1}{c|}{}& & & & & &

\\
& \multicolumn{1}{c|}{}& & & &  & & \multicolumn{1}{c|}{}& & & & & &
\\
\hline
\end{tabular}
\end{center}
\caption{The condition number of $\mathbf{B}_{E}^{type} \mathbf{A}$ and the number of iterations of the preconditioned conjugate gradient method (in parentheses) for $type=\RN{2}$ with different values for $\alpha_b$, $\alpha_c$, and $\alpha_i$. For $6 \times 6$ and $9 \times 9$ subdomains, $h$ and $h^*$ belong to $\{ 1/36,1/54 \}$ and $\{ 1/54,1/81 \}$, respectively. In addition, selecting the number of eigenfunctions for each subdomain to construct the enrichment coarse space is based on the adaptive enrichment, where the threshold is $50$.}
\label{Table1}
\end{table}

In Table~\ref{Table1}, we set different values for $\alpha_b$, $\alpha_c$, and $\alpha_i$ for different number of subdomains, for instance, see Figure~\ref{first alpha} including $6\times6$ subdomains and the coarse mortar case. We implement the additive average Schwarz method with the adaptive enrichment coarse space $type=\RN{2}$, where the threshold is 50, i.e., all eigenfunctions associated with all eigenvalues larger than $50$ are considered. Note that the condition number of the non-enrichment preconditioner is very large. For instance, in the case  $6\times 6$ subdomains and coarse mortar, it is $1.32\text{e}7$. We first observe that the condition numbers of the enrichment preconditioners are proportional to the ratio $H/h$ and independent of the number of subdomains. We also observe that if we consider the same number of eigenfunctions for each subdomain based on a proper threshold,  consequently, the condition number estimates are independent of values of $\alpha_b$, $\alpha_c$ and $\alpha_i$.

Table~\ref{Table2} demonstrates only the condition number estimates and iteration numbers of the preconditioned conjugate gradient method (in parentheses) in the cases $type=\RN{2}$, coarse mortar, $6\times 6$ and $9\times9$ subdomains, where $\alpha_b=1$, $\alpha_c=1\text{e}4$, and $\alpha_i=1\text{e}6$. For the other cases and also different values of jumps in the coefficient $\alpha(x)$, we have similar results. Also, the construction of the enrichment coarse space is relied on the fixed number of eigenfunctions for each subdomain varying from $0$ to $7$.

\begin{table}[!ht]
\begin{center}
\begin{tabular}{|c|cccc c c c c c| }
\cline{1-10}
\multicolumn{2}{|c|}{\multirow{2}{*}{ADD}} & \multicolumn{8}{c|}{\tiny  Number of Eigenfunctions in Each Subdomain} \\
\cline{3-10}
\multicolumn{2}{|c|}{} & 0 & 1 & 2 & 3 & 4 & 5 & 6 & 7\\
\hline
\multirow{6}{*}{\begin{sideways} {\tiny Number of Subdomains \text{ }} \end{sideways}}  & \multicolumn{1}{c|}{ \multirow{4}{*}{$6 \times 6$} } & \multirow{3}{*}{$1.32\text{e}7$}&\multirow{3}{*}{$3.25\text{e}6$}& \multirow{3}{*}{$2.65\text{e}5$}& \multirow{3}{*}{$4.10\text{e}4$}& \multirow{3}{*}{$4.12\text{e}3$} & \multirow{3}{*}{$5.7\text{e}1$} & \multirow{3}{*}{$5.6\text{e}1$} & \multirow{3}{*}{$5.58\text{e}1$}\\ 
&\multicolumn{1}{c|}{} & & & & & & & &\\
& \multicolumn{1}{c|}{}& \multirow{1}{*}{$(1756)$}  & \multirow{1}{*}{$(1168)$}   & \multirow{1}{*}{$(584)$}  & \multirow{1}{*}{$(295)$} &  \multirow{1}{*}{$(112)$}  &  \multirow{1}{*}{$(50)$} &  \multirow{1}{*}{$(48)$}  & \multirow{1}{*}{$(47)$}\\
& \multicolumn{1}{c|}{\multirow{6}{*}{$9\times9$}} & \multirow{5}{*}{$\text{2.29e7}$}& \multirow{5}{*}{$\text{7.54e6}$}&\multirow{5}{*}{$1.41\text{e}6$}& \multirow{5}{*}{$3.82\text{e}5$}& \multirow{5}{*}{$1.23\text{e}5$}&\multirow{5}{*}{$6.21\text{e}3$}&\multirow{5}{*}{$5.80\text{e}2$} &\multirow{5}{*}{$7.38\text{e}1$}\\
\cline{2-10} 
& \multicolumn{1}{c|}{}& \multirow{5}{*}{$(4981)$}  & \multirow{5}{*}{$(3889)$}   & \multirow{5}{*}{$(2368)$}  & \multirow{5}{*}{$(1180)$} &  \multirow{5}{*}{$(353)$}  &  \multirow{5}{*}{$(118)$} &  \multirow{5}{*}{$(63)$} & \multirow{5}{*}{$(50)$} \\
& \multicolumn{1}{c|}{}& & & & & & & &\\
& \multicolumn{1}{c|}{}& & & & & & & &

\\
& \multicolumn{1}{c|}{}& & & & & & & &
\\
\hline
\end{tabular}
\end{center}
\caption{The implementation of the additive average Schwarz method, $type=\RN{2}$, with the  fixed number of eigenfunctions for each subdomain to estimate the condition number of $\mathbf{B}_{E}^{type} \mathbf{A}$ and the number of iterations of the preconditioned conjugate gradient method (in parentheses) for $6 \times 6$ ($h=1/36$, $h^*=1/54$) and $9 \times 9$ ($h=1/54$, $h^*=1/81$) subdomains and the coarse mortar case. Further, the  distribution of jumps in $\alpha(x)$ are $\alpha_{b}=1$, $\alpha_c=1\text{e}4$ and $\alpha_i=1\text{e}6$.}
\label{Table2}
\end{table}

\begin{table}[!ht]
\begin{center}
\begin{tabular}{|c|cc|c|cc c c c cccccccc | }
\cline{1-16}
\multicolumn{4}{|c|}{\multirow{3}{*}{ADD}} & \multicolumn{12}{c|}{ \tiny{ Different Values for} $\alpha_b$, $\alpha_c$ and $\alpha_i$ }\\
\cline{5-16}
\multicolumn{4}{|c|}{} &\multicolumn{6}{|c|}{ $1$, $1\text{e}3$,  $1\text{e}4$} & \multicolumn{6}{|c|}{$1$,  $1\text{e}4$,  $1\text{e}6$ } \\
\cline{5-16}
\multicolumn{4}{|c|}{} &\multicolumn{3}{|c|}{\tiny Coarse Mortar }& \multicolumn{3}{|c|}{\tiny Fine Mortar } & \multicolumn{3}{|c|}{\tiny Coarse Mortar}  &  \multicolumn{3}{|c|}{\tiny Fine Mortar }\\
\hline
\multirow{6}{*}{\begin{sideways} {\tiny Number of Subdomains \text{ }} \end{sideways}}  & \multicolumn{1}{c|}{ \multirow{4}{*}{$6 \times 6$} }& \multirow{4}{*}{\begin{sideways} {\tiny $type$} \end{sideways}} & \multicolumn{1}{c|}{ \multirow{3}{*}{ $\RN{1}$}} & \multicolumn{3}{c}{ \multirow{3}{*}{$625$ }}&\multicolumn{3}{c|}{ \multirow{3}{*}{ $525$}}& \multicolumn{3}{c}{\multirow{3}{*}{$635$ }} &  \multicolumn{3}{c|}{\multirow{3}{*}{ $534$}} \\ 
&\multicolumn{1}{c|}{} & \multicolumn{1}{|c|}{}&\multicolumn{1}{c|}{} & & & & & & \multicolumn{1}{c|}{}& & & & & & \multicolumn{1}{c|}{}\\
& \multicolumn{1}{c|}{}& & \multicolumn{1}{c|}{\multirow{1}{*}{$\RN{2}$}} & \multicolumn{3}{c}{\multirow{1}{*}{$75$}}    & \multicolumn{3}{c|}{\multirow{1}{*}{$80$}} & \multicolumn{3}{c}{\multirow{1}{*}{$75$}} & \multicolumn{3}{c|}{\multirow{1}{*}{$80$}}   \\
& \multicolumn{1}{c|}{\multirow{6}{*}{$9\times9$}} & \multirow{6}{*}{\begin{sideways} {\tiny $type$} \end{sideways}} &\multicolumn{1}{c|}{\multirow{5}{*}{$\RN{1}$ }} & \multicolumn{3}{c}{\multirow{5}{*}{$1988$ }}& \multicolumn{3}{c|}{\multirow{5}{*}{$1680$}} & \multicolumn{3}{c}{\multirow{5}{*}{$1992$ }}  &  \multicolumn{3}{c|}{\multirow{5}{*}{$1681$}}\\
\cline{2-16} 
& \multicolumn{1}{c|}{} & & \multicolumn{1}{c|}{\multirow{5}{*}{$\RN{2}$}} & \multicolumn{3}{c}{\multirow{5}{*}{$240$}}    &   \multicolumn{3}{c|}{\multirow{5}{*}{$233$}}  &   \multicolumn{3}{c}{\multirow{5}{*}{$240$}}   &    \multicolumn{3}{c|}{\multirow{5}{*}{$233$}}\\
& \multicolumn{1}{c|}{}& & & & & & & & \multicolumn{1}{c|}{}& & & & & & \multicolumn{1}{c|}{}\\
& \multicolumn{1}{c|}{}& & & & & & & & \multicolumn{1}{c|}{}& & & & & & \multicolumn{1}{c|}{}

\\
& \multicolumn{1}{c|}{}& & & &  & & & &\multicolumn{1}{c|}{} & & & & & &\multicolumn{1}{c|}{}
\\
\hline 
\end{tabular}
\end{center}
\caption{The total numbers of the eigenfunctions associated with the eigenvalues greater than $50$ to enrich the coarse spaces  used in the additive average Schwarz method for both types $\RN{1}$ and $\RN{2}$.}
\label{Table3}
\end{table}

It is conclusive that both approaches to enrich the standard coarse space, i.e., using the given threshold and imposing the fixed numbers of eigenfunctions for each subdomain, lead to similar results. This fact can be viewed by comparing, for instance, the third column of Table~\ref{Table1} with the last column of Table~\ref{Table2}. 

\begin{figure}[!ht]
\centering
\setlength{\unitlength}{0.1\textwidth}
\vspace{9.5cm}
\begin{picture}(15,6.5)

\put(-.6,7.5){ \includegraphics[width=6.4cm]{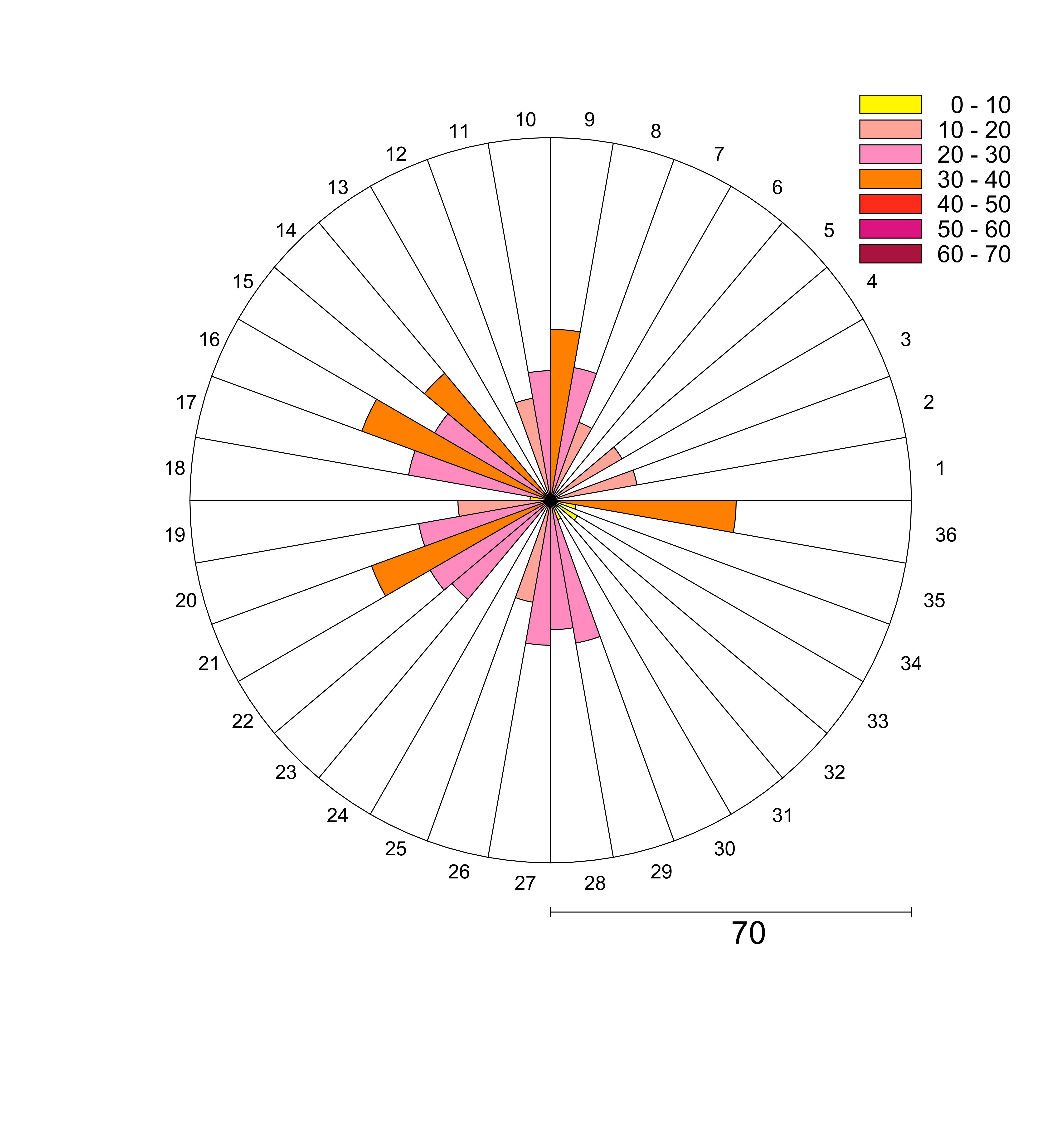}}
  \put(.6,8.1){$(a)$ Fine mortar and $type=\RN{1}.$}
  \put(4.7,7.15){ \includegraphics[width=7.15cm]{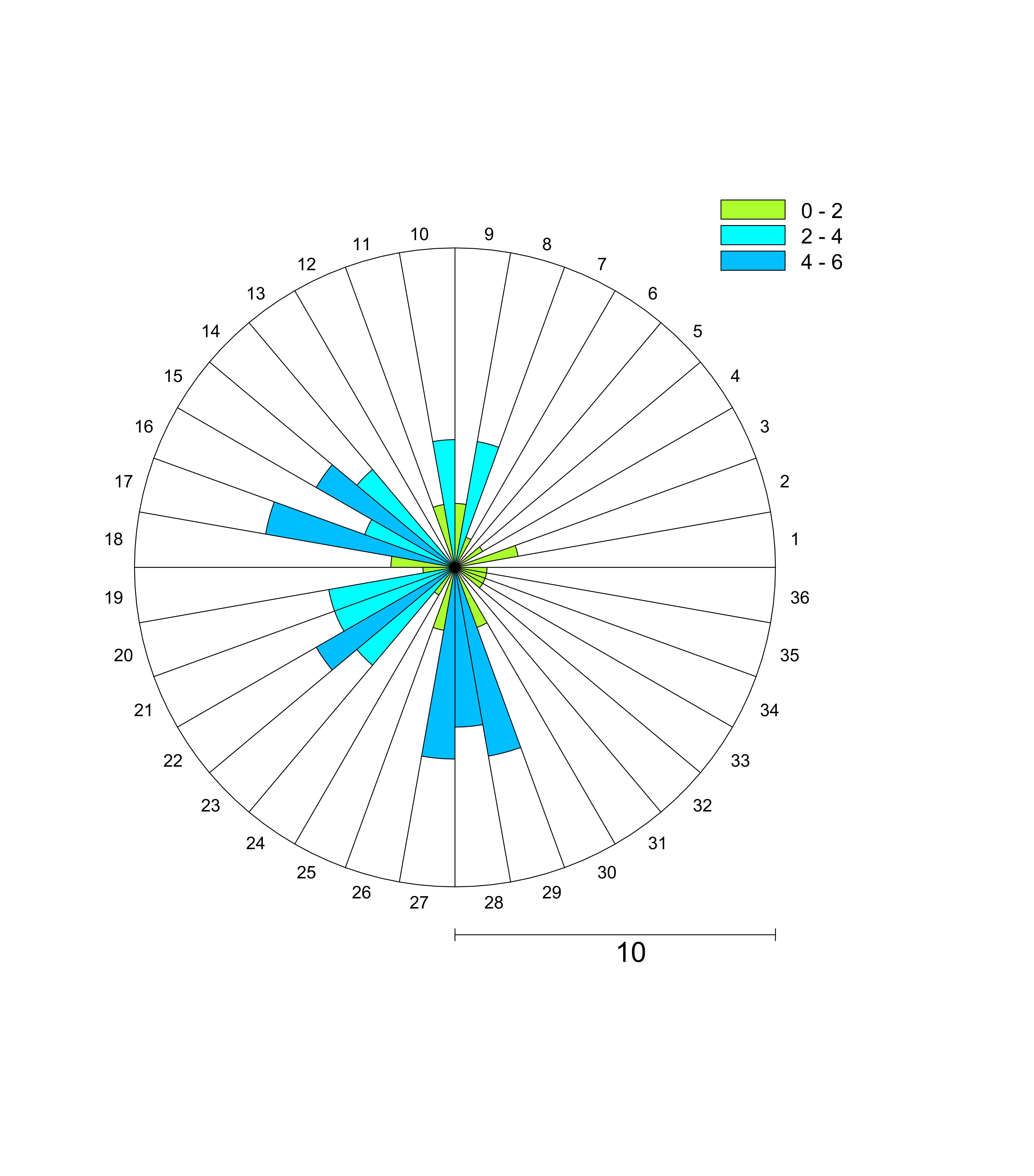}}
  
\put(-.6,1.5){ \includegraphics[width=6.5cm]{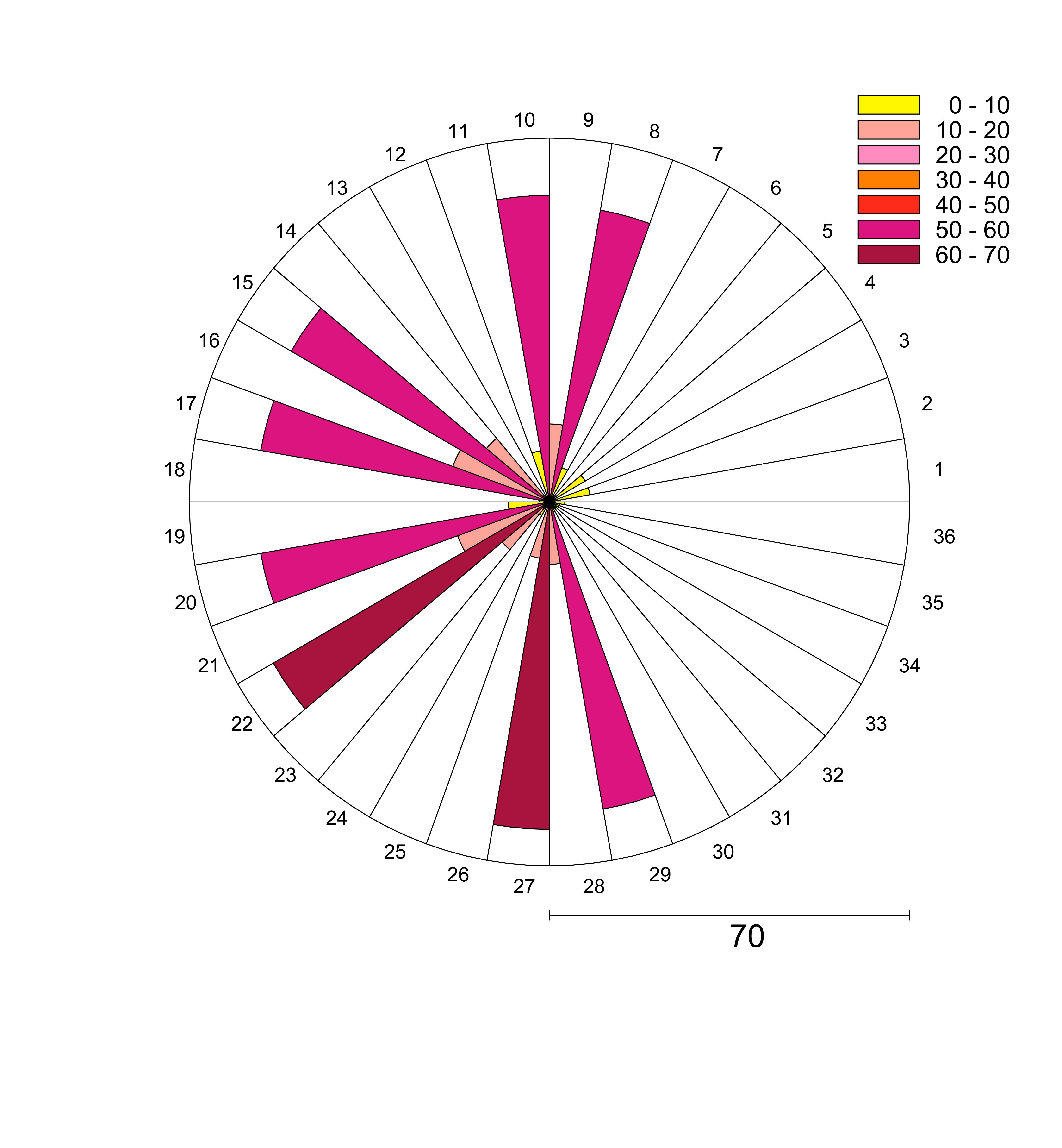}}
  \put(0.5,2.1){$(c)$ Coarse mortar and $type=\RN{1}.$}
  \put(4.7,1.25){ \includegraphics[width=7.2cm]{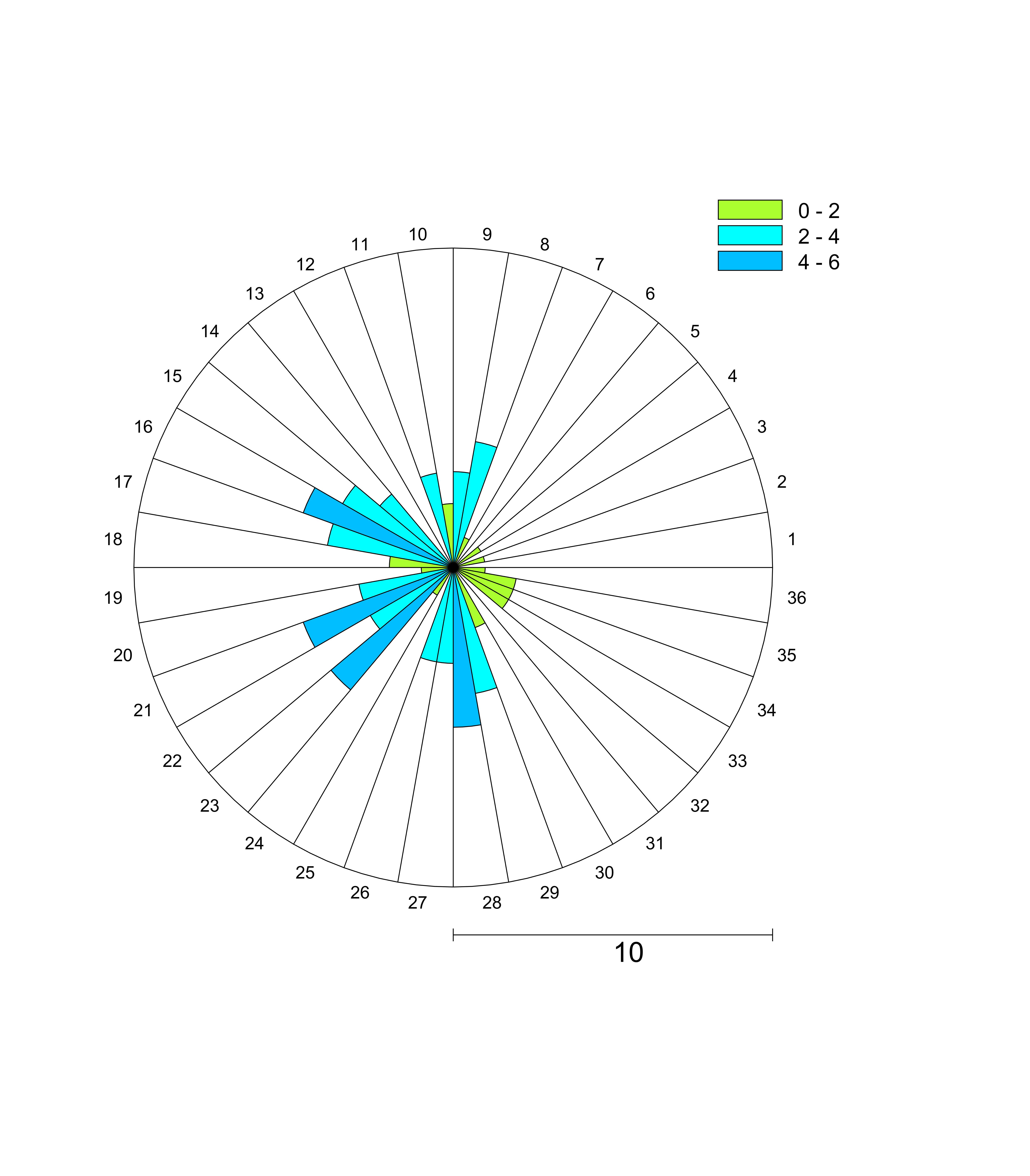}}
  \put(5.6,2.1){$(d)$ Coarse mortar and $type=\RN{2}.$}
  
  \put(5.7,8.){$(b)$ Fine mortar and $type=\RN{2}.$} 
  
  \end{picture}
  \vspace{-2.5cm}
  \centering
  \caption{The histograms of the number of eigenfunctions associated with  the  eigenvalues greater than $50$ corresponding to the partition of $\Omega=[0,1]\times[0,1]$ into $6\times 6$ subdomains, where $\alpha_b=1$, $\alpha_c=1\text{e}4$ and $\alpha_i=1\text{e}6$. The number of subdomains are ordered \!anticlockwise around the polar histograms from 1 to 36. For $type=\RN{1}$ and $\RN{2}$, the number of eigenfunctions are grouped into the ranges $0$ to $70$ and $0$ to $10$, respectively. The largest numbers of the eigenfunctions in $(a)-(d)$ are $39$, $6$, $63$ and $5$, respectively.}
  \label{histogram type I II}
\end{figure}

Table~\ref{Table3} gives the total number of required  eigenfunctions for the adaptive enrichment coarse spaces in the cases $6\times 6$ and $9 \times 9$ subdomains, $type=\RN{1}$ and $\RN{2}$ and different values for $\alpha(x)$. As we can see from this table, solving the second type of the generalized eigenvalue problem \eqref{GenEigPro} is more efficient than solving the first type. To analyze the distribution of eigenfunctions,  Figure~\ref{histogram type I II} contains the polar histograms for both coarse and fine mortars in the case $6\times 6$ subdomains, and $type \in \{\RN{1}, \RN{2}\}$. It clearly shows that we need only consider a few eigenfunctions for each subdomain in the case $type=\RN{2}$ compared to $type=\RN{1}$.

\section{Conclusion}
In this paper, we have employed the additive average Schwarz method with the enrichment coarse spaces to solve second order elliptic boundary problem coupled  with very large jumps in the function $\alpha(x)$, where the finite element discretization of that problem has been based on the nonmatching triangulations across the subdomains interfaces. We have proved that the condition number estimates for the produced preconditioners by the additive average Schwarz methods for $type \in \{ \RN{1},\RN{2} \}$ are proportional to the ratio $H/h$ and independent of the number of subdomains. Besides, we have compared the numerical  results to conclude that, in practice, $type=\RN{2}$ has much better performance than $type=\RN{1}$.

\section*{Acknowledgements}
The work of Leszek Marcinkowski  was partially supported by Polish Scientific Grant: National Science Center:  2016/21/B/ST1/00350.


\begin{thebibliography}{}
%
%





\bibitem{Arbogast-2000} 
T. Arbogast, L. C. Cowsar, M. F. Wheeler and I. Yotov, Mixed finite element methods on nonmatching multiblock grids, 
{\it SIAM J. Numer. Anal.} 37(4), 1295--1315, (2000)

\bibitem{Arbogast-2007}
T. Arbogast, G. Pencheva, M. F. Wheeler and I. Yotov, A multiscale mortar mixed finite element method,
 {\it Multiscale Model. Simul.} 6(1), 319--346, (2007)

\bibitem{Arbogast-2013}
T. Arbogast and H. Xiao, A multiscale mortar mixed space based on homogenization for heterogeneous elliptic problems, 
{\it SIAM J. Numer. Anal.} 51(1), 377--399, (2013)
\bibitem{Bernardi-2005} 
C. Bernardi, Y. Maday and F. Rapetti, Basics and some application of mortar element method,
GAMM-Mitt, 28, 97--123, (2005)
\bibitem{Brenner-1996}
S. C. Brenner, 
Two-level additive Schwarz preconditioners for nonconforming finite element methods,
{\it Math. Comp.} 65, 897--921, (1996)
\bibitem{Brezina-1999}
 M. Brezina, C. Heberton, J. Mandel and P. Van\v{e}k, An iterative method with convergence rate
chosen a priori, Tech. Report 140, Center for Computational Mathematics, University of Colorado 
Denver, Denver, (1999)
\bibitem{Bjorstad-1996}
P. E. Bj{\o}rstad, M. Dryja and E. Vainikko, Parallel implementation of a Schwarz domain decomposition algorithms, in Applied Parallel Computing in Industrial Problems and Optimization (J. Wasniewski, J. Dongara, K. Madsen, and D. Olsem, eds.) Lecture Notes in Computer Science, Springer, 1184, 141--157, (1996)
\bibitem{Bjorstad-2003}
P. E. Bj{\o}rstad, M. Dryja and T. Rahman, Additive Schwarz methods for elliptic mortar finite element
problems, 
{\it Numer. Math.} 95, 427--457, (2003)
\bibitem{Chan-1994}
T. F. Chan and T. P. Mathew, Domain decomposition algorithms, 
{\it Acta Numerica.} 3, 61--143, (1994)
\bibitem{Ciarlet-1978}
P. G. Ciarlet,
The Finite Element Method for Elliptic Problems,
North-Holland, Amsterdam, (1978)
\bibitem{Cowsar-1993}
 L. C. Cowsar, Domain decomposition methods for nonconforming finite elements spaces of
 Lagrange-type, Proceedings of the Sixth Copper Mountain Conference on Multigrid Methods,
 NASA Conference Publication 3224, 93--109, (1993)
 \bibitem{Dryja-1987}
M. Dryja and O. B. Widlund, 
An additive variant of the Schwarz alternating method for the case of many subregions, 
Technical Report 339, also Ultracomputer Note 131, Department of Computer Science, Courant Institute, New York University, (1987)
\bibitem{Dryja-2004}
M. Dryja, A. Gantner, O. B. Widlund and B. I. Wohlmuth, Multilevel
additive Schwarz preconditioner for nonconforming mortar finite element methods, {\it J. Numer. Math.}  12, 23--38, (2004)
\bibitem{Dryja-2010}
M. Dryja and M. Sarkis, Additive average Schwarz methods for discretization of elliptic problems with highly discontinuous coefficients,
{\it Comput. Methods Appl. Math}. 10, 164--176, (2010)
\bibitem{Efandiev-2011}
Y. Efendiev, J. Galvis and P. S. Vassilevski, Spectral element agglomerate algebraic multigrid methods for
elliptic problems with high contrast coefficients. In: Y. Huang, R. Kornhuber, O. B. Widlund, J. Xu (eds.)
Domain Decomposition Methods in Science and Engineering XIX, Lecture Notes in Computational
Science and Engineering, vol. 78, pp. 407–414. Springer, Berlin (2011)
\bibitem{Eikeland-2019} 
E. Eikeland, L. Marcinkowski and T. Rahman, Overlapping Schwarz methods with adaptive coarse spaces for multiscale problems in 3D, 
{\it Numer. Math.} 142, 103--128, (2019)
\bibitem{Feng-2002}
X. Feng and T. Rahman, An additive average Schwarz method for the plate bending problem,
{\it J. Numer. Math.} 10, 109--125, (2002) 
\bibitem{Gander-2015} 
M. J. Gander, A. Loneland and T. Rahman, Analysis of a new harmonically enriched multiscale coarse space for domain decomposition methods, arXiv preprint arXiv:1512.05285 (2015)

\bibitem{Heinlein-2019}
A. Heinlein, A. Klawonn, J. Knepper and O. Rheinbach, 
Adaptive GDSW coarse spaces for overlapping Schwarz methods in three dimensions,
{\it SIAM J. Sci. Comput.} 41(5), A3045--A3072, (2019)
\bibitem{Kim-2015} 
H. H. Kim and E. T. Chung, A BDDC algorithm with enriched coarse spaces for two-dimensional elliptic problems with oscillatory and high contrast coefficients,
{\it Multiscale Model. Simul.} 13, 571--593, (2015)
\bibitem{Kim-2017}
H. H. Kim, E. Chung and J. Wang, 
BDDC and FETI-DP preconditioners with adaptive coarse spaces for three-dimensional elliptic problems with oscillatory and high contrast coefficients,
{\it J. Comput. Phys.} 349, 191--214, (2017)
\bibitem{Klawonn-2014} 
A. Klawonn, M. Lanser, P. Radtke and O. Rheinbach, On an adaptive coarse space and on nonlinear domain decomposition. In: J. Erhel,
M. J. Gander, L. Halpern, G. Pichot, T. Sassi, O. B. Widlund (eds.) Domain Decomposition Methods in Science and Engineering XXI,  Lecture Notes in Computational Science and Engineering, vol. 98, pp. 71-83. Springer, Switzerland (2014)
\bibitem{Klawonn-2015}
A. Klawonn, P. Radtke and O. Rheinbach, FETI-DP methods with an adaptive coarse space,
{\it SIAM J. Numer. Anal.} 53(1), 297--320, (2015)
\bibitem{Loneland-2016}
A. Loneland, L. Marcinkowski and T. Rahman,
 Additive average Schwarz method for a Crouzeix-Raviart finite volume element discretization of elliptic
problems with heterogeneous coefficients, 
{\it Numer. Math.} 134, 91--118, (2016)
\bibitem{Mandel-2012} 
J. Mandel, B. Sousedk and J. \v{S}\'{i}stek, Adaptive BDDC in three dimensions,
{\it Math. Comput. Simulation.} 82, 1812--1831, (2012)
\bibitem{Marcinkowski-2018}
 L. Marcinkowski and T. Rahman, Additive average Schwarz with adaptive coarse spaces: scalable algorithms for multiscale problems, 
 {\it Electron. Trans. Numer. Anal.} 49, 28--40., (2018)
\bibitem{Matsokin-1985}
A. M. Matsokin and S. V. Nepomnyaschikh, 
A Schwarz alternating method in a subspace.
{\it Sov. Math.} 29(10), 78--84, (1985)
\bibitem{Rahman-2005}
T. Rahman, X. Xu and R. Hoppe, Additive Schwarz methods for the Crouzeix-Raviart mortar finite element for elliptic problems with discontinuous coefficients, 
{\it Numer. Math.} 101(3), 551--572, (2005)
\bibitem{Sarkis-1993}
M. Sarkis, 
Two-level Schwarz methods for nonconforming finite elements and discontinuous coefficients.
Proceedings of the Sixth Copper Mountain Conference on Multigrid Methods, NASA Conference Publication 3224, 543--565, (1993)
\bibitem{Sarkis-1997}
M. Sarkis, Nonstandard coarse spaces and Schwarz methods for elliptic problems with discontinuous
coefficients using non-conforming elements, 
{\it Numer. Math.} 77(3), 383--406, (1997)
\bibitem{Seshiyer-2000} 
P. Seshaiyer and M. Suri, Uniform $hp$ convergence results for the mortar finite element method,
{\it Math. Comp.} 69, 521--546, (2000)
\bibitem{Smith-book-1996}
B. F. Smith, P. E. Bjørstad and W. D. Gropp, Domain Decomposition: Parallel Multilevel Methods for Elliptic Partial Differential Equations, Cambridge University Press, (1996).
\bibitem{Spillane-2013}
N. Spillane and D. J. Rixen, Automatic spectral coarse spaces for robust finite element tearing and interconnecting and balanced domain
decomposition algorithms,
{\it Int. J. Numer. Methods. Eng.} 95(11), 953--990, (2013)
\bibitem{Spillane-2014}
N. Spillane, V. Dolean, P. Hauret, F. Nataf, C. Pechstein and R. Scheichl, Abstract robust coarse spaces for systems of PDEs via generalized eigenproblems in
the overlaps,
{\it Numer. Math.} 126(4), 741--770, (2014)
\bibitem{Tallec-1991}
P. Le Tallec, Y. H. De Roeck and M. Vidrascu, Domain decomposition methods for large linearly elliptic three-dimensional problems,
{\it J. Comput. Appl. Math.} 34, 93--117, (1991)
\bibitem{Toselli-2005}
A. Toselli and O. B. Widlund, Domain decomposition methods-algorithms and theory, volume 34 of Springer Series in Computational Mathematics, Springer-Verlag, Berlin, (2005)
\bibitem{Xu-1998}
J. Xu and J. Zou, Some nonoverlapping domain decomposition methods,
{\it SIAM Rev.} 40, 857--914, (1998)


\end{thebibliography}


\end{document}